\documentclass[review,hidelinks,onefignum,onetabnum]{siamart251216}

\usepackage{amssymb}
\usepackage{lineno}

\usepackage{amsmath,mathtools,bm}
\usepackage{nicefrac}
\usepackage{subfig}
\usepackage{adjustbox}   
\usepackage{tikz}
\usetikzlibrary{arrows,shapes,backgrounds}

\usepackage{pgfplots}
\usetikzlibrary{arrows}
\pgfplotsset{compat=1.10}
\usetikzlibrary{calc}    

\definecolor{MPIblue}{RGB}{51,165,195}
\definecolor{MPIgrey}{RGB}{135,135,141}
\definecolor{MPIgreen}{RGB}{0,118,117}
\definecolor{MPIred}{RGB}{120,0,75}
\definecolor{MPIsand}{RGB}{236,233,212}
    
\newtheorem{prop}[theorem]{Proposition}

\newcommand{\coloneqq}{\mathrel{\mathop:}=}
\newcommand{\D}[1]{\, \mathrm{d} #1}

\newcommand*{\Shat}{\skew{3}{\hat}{\bm{S}}}

\newcommand\Tstrut{\rule{0pt}{2.6ex}}         
\newcommand\Bstrut{\rule[-0.9ex]{0pt}{0pt}}   
\newcommand\BBstrut{\rule[-0.9ex]{0pt}{0pt}}   

\headers{Preconditioning for Diffuse Interface Tumor Growth}{J. Bosch, P. \c{C}\.{i}lo\u{g}lu, C. Greif, and M. Stoll}

\title{Preconditioning for Diffuse Interface Tumor Growth Models}

\author{
Jessica Bosch\thanks{Department of Computer Science, The University of British Columbia, Vancouver, BC, V6T 1Z4, Canada (\email{jbosch@cs.ubc.ca},  \email{greif@cs.ubc.ca}).} 
\and Pel\.{i}n \c{C}\.{i}lo\u{g}lu\thanks{Faculty of Mathematics, Chemnitz University of Technology, 09107, Chemnitz, Germany   (\email{pelin.ciloglu@mathematik.tu-chemnitz.de}, \email{martin.stoll@mathematik.tu-chemnitz.de}).}
\and Chen Greif\footnotemark[1]
\and Martin Stoll\footnotemark[2]
}

\begin{document}
\nolinenumbers
\maketitle

\begin{abstract}
Computational models of tumor growth form a contemporary area of study. Diffuse interface models have been proposed as an important tool for cancer modeling. In this work, we consider preconditioners for numerical simulations of two tumor growth models based on diffuse interface models. For numerical efficiency, we propose block preconditioners that rely on using effective Schur-complement approximations. We show that these methods combined with state-of-the-art adaptive finite element discretizations lead to robust simulations in two and three dimensions. In extensive numerical experiments, we show that the proposed methods show robust convergence behavior.
\end{abstract}

\begin{keyword}
Cahn--Hilliard equation, diffuse interface model, tumor growth, preconditioning, Schur complement approximation, finite element discretizations
\end{keyword}

\begin{MSCcodes}
35K55, 35Q92, 65F08, 65F10, 65M60
\end{MSCcodes}

\section{Introduction}
\label{intro} 
\sloppy
Understanding tumor growth is a challenging multiscale problem involving complex interactions between biological, chemical, and mechanical processes. Mathematical models provide a valuable tool for studying these processes, enabling researchers to investigate tumor evolution and explore treatment strategies in silico before conducting costly, possibly dangerous, and time-consuming laboratory or animal experiments~\cite{HawZT12}. In the long term, such models may also contribute to personalized medicine by helping predict the qualitative progression of individual tumors.

Among the many mathematical approaches to tumor growth, diffuse interface models have been proposed as a tool with an important significance to cancer modeling given their flexibility and ability to easily incorporate topological changes.  In these models, tumor and healthy tissue are represented as two phases separated by a thin but diffuse interface described by a phase-field variable. This avoids explicit interface tracking while naturally accommodating the mentioned topological changes. The governing equations are typically of Cahn--Hilliard type, derived as gradient flows of Ginzburg--Landau energy functionals~\cite{CahnHilliard1958}, and can incorporate important biological mechanisms such as nutrient-driven proliferation, apoptosis, chemotaxis, and active transport. Early contributions by Cristini et al.~\cite{cristini2009nonlinear, CristiniLowengrub2010} and Wise et al.~\cite{Wise2008} considered multiphase tumor growth with nutrient transport and Darcy flow. Rigorous mathematical analysis of Cahn--Hilliard tumor growth models was later carried out in~\cite{ColliGilardiHilhorst2015, FrigeriGrasselliRocca2015, HikKNZ15}. Garcke et al.~\cite{GarLSS16} proposed a thermodynamically consistent Cahn--Hilliard--Darcy model incorporating chemotaxis and active transport. Further analytical and numerical developments can be found in~\cite{garcke2022viscoelastic,GarckeLam2017dirichlet,GarckeLam2017wp,GarckeLamNurnbergSitka2018,GarckeLamSignori2021,garcke2022numerical}.

In this paper, we consider two diffuse interface tumor growth models. The first is the model of Hilhorst--Kampmann--Nguyen--Zee (HKNZ)~\cite{HikKNZ15}, which couples a Cahn--Hilliard equation with a reaction--diffusion equation for the nutrient. The second is the Garcke--Lam--Sitka--Styles (GLSS) model~\cite{GarLSS16} in both its smooth and nonsmooth forms, which additionally accounts for chemotaxis and active transport.

The main computational challenge in diffuse interface tumor growth modeling is the efficient solution of the large and sparse linear systems arising at each time step. For the smooth GLSS model, Newton method produces a sequence of nonsymmetric saddle--point systems, while a Moreau--Yosida regularization combined with a semi-smooth Newton method is employed. Efficient preconditioners for Cahn--Hilliard systems have been developed~\cite{ BKSW14, BosS15a, BSB14, boyanova2012efficient,boyanova2014efficient},  and, more recently, for Cahn--Hilliard--Navier--Stokes systems ~\cite{BKS17, PCiloglu_2025, kay2007efficient}.

Building on these ideas, we develop block preconditioners based on Schur complement approximations for both the smooth and nonsmooth Cahn--Hilliard tumor growth  GLSS and HKNZ models. To the best of our knowledge, dedicated preconditioning strategies for Cahn--Hilliard-based tumor growth models have not been reported in the existing literature. This work aims to fill this gap by introducing efficient and robust block preconditioners for these models. Through extensive two- and three-dimensional numerical experiments, we demonstrate that the proposed preconditioners are robust with respect to model parameters, time-step sizes, and mesh refinement, making them suitable for the efficient solution of large-scale tumor growth simulations.

The paper is organized as follows. Section~\ref{sec:CHtumor} introduces the mathematical models of tumor growth. The time discretization is given in Section~\ref{sec:time} and then the semi-smooth Newton method is described in Section~\ref{sec:newton}. Section~\ref{sec:FEM} presents the finite element discretization and the resulting linear systems. The proposed block preconditioners for the smooth and nonsmooth models are developed in Sections~\ref{sec:prec}. Numerical results are presented in Section~\ref{sec:results}, followed by concluding remarks in Section~\ref{sec:concl}.

\section{Diffuse interface tumor growth models}
\label{sec:CHtumor}
Let $\Omega\subset\mathbb{R}^{d}$, where $d\in\{2,3\}$, be a bounded domain that contains a two-component mixture consisting of tumor and healthy cells. We denote the tumor cell phase by $u\colon\Omega\to\mathbb{R}$. A tumorous phase is characterized by $u\approx 1$ and  a healthy cell phase by $u\approx -1$. The thickness of the interface between those two phases is proportional to the parameter $\varepsilon>0$.

We denote the nutrient phase by $\sigma\colon\Omega\to\mathbb{R}$. A nutrient-rich extracellular water phase is characterized by $\sigma\geq 0$.  A strong separation between nutrient-rich and nutrient-poor phases is not anticipated. Hence there is neither a surface energy nor a double well potential with respect to the nutrient phase included in the energy functional. The dynamics of the nutrient species is, without interaction, simply governed by diffusion.

\subsection{Hilhorst--Kampmann--Nguyen--Zee (HKNZ) model}
\label{sec:Hilhorst-model}
In Hilhorst et al.~\cite{HikKNZ15}, the unknown pair $(u,\sigma)$ is a dissipative gradient flow for the energy functional
\begin{align*}
 \mathcal{E}_{\text{HKNZ}}(u,\sigma)=\int_{\Omega}{\frac{\varepsilon}{2}|\nabla{u}|^{2}+\frac{1}{\varepsilon}\psi(u)+\frac{1}{2s}\sigma^2\,\D{\bf{x}}}\,,
\end{align*}
where $\psi(u)=0.25\,(1-u^2)^2$ is a potential function with two equal minima at $\pm 1$  and  $s>0$ a small regularization parameter.

The diffuse interface tumor growth model in Hilhorst et al.~\cite{HikKNZ15}, which is in the following denoted by the HKNZ model, is given by
\begin{subequations}   
\begin{align} 
  u_t &= \Delta{w}+\varepsilon^{-1}p(u)(\sigma-s w)\,,\label{HikKN1}\\
  w &= -\varepsilon\Delta{u}+\varepsilon^{-1}\psi'(u)\,,\label{HikKN2}\\
  \sigma_t &= \Delta{\sigma}-\varepsilon^{-1}p(u)(\sigma-s w)\,,\label{HikKN3}\\
  \nabla{u}\cdot\mathbf{n}&=\nabla{\sigma}\cdot\mathbf{n}=\nabla{w}\cdot\mathbf{n}=0\quad\textup{on}\ \partial\Omega\,,\label{HikKN4}
\end{align}
\end{subequations}
where $w\colon\Omega\to\mathbb{R}$ is the chemical potential. The growth function $p(u)$ is defined by
\begin{align}\label{eqn:p}
  p(u):=\left\{\begin{array}{rl}p_0(1-u^2) &\  u\in[-1,1]\,,\\0 &\ \textup{elsewhere},\end{array}\right.
\end{align}
where $p_0>0$ is a proliferation growth parameter. 

The second terms on the right-hand side of (\ref{HikKN1}) and (\ref{HikKN3}) are the mass exchange terms, which are motivated by the linear phenomenological constitutive laws for chemical reactions. It is anticipated that there will be only growth when there is some level of nutrient available. In addition, the choice of $p(u)$ in \eqref{eqn:p} helps prevent growth if, due to numerical discretization, the solution takes on absolute values larger than $1$.

In the next section, we consider another recently proposed model which in addition introduces terms that reflect chemotaxis, which is the active movement of the tumor colony towards nutrient source.

\subsection{Garcke--Lam--Sitka--Styles (GLSS) model}
\label{sec:Garcke-model}
In Garcke et al.~\cite{GarLSS16}, the unknown pair $(u,\sigma)$ is a dissipative gradient flow for the energy functional
\begin{align*}
 \mathcal{E}_{\text{GLSS}}(u,\sigma)=\int_{\Omega}{\frac{\beta\varepsilon}{2}|\nabla{u}|^{2}+\frac{\beta}{\varepsilon}\psi(u)+\frac{\chi_u}{2\lambda}\sigma^2+\chi_u \sigma(1-u)\,\D{\bf{x}}}.
\end{align*}

As before, $\psi(u)$ is a potential function with two equal minima at $\pm 1$. We consider either the smooth potential $\psi(u)=0.25\,(1-u^2)^2$ or the nonsmooth double-obstacle potential
\begin{align*}
\psi(u)=\left\{\begin{array}{rl}\psi_0(u)=\frac{1}{2}(1-u^{2}) & |u|\leq 1,\\ \infty & |u|>1.\end{array}\right.
\end{align*}

The positive constant $\beta$ denotes the surface tension and $\chi_u \geq 0$ can be seen as a transport mechanism such as chemotaxis and active uptake. With $\lambda>0$, we can switch off effects of active transport by letting $\lambda\to 0$, while preserving the effects of chemotaxis. Such a term is not present in the HKNZ model. 

The diffuse interface tumor growth model in Garcke et al.~\cite{GarLSS16}, which we call the GLSS model, is given by
\begin{subequations}
\begin{align} 
  u_t &= \text{div}\left(0.5(1+u)^2\nabla{w}\right)+\mathcal{P}\sigma(u+1)-\mathcal{A}(u+1)\,,\label{GarLSS161}\\
  w &= -\beta\varepsilon\Delta{u}+\beta\varepsilon^{-1}\psi'(u)-\chi_{u}\sigma\,,\label{GarLSS162}\\
  \sigma_t &= \text{div}\left(\mathcal{D}(u)\nabla{\sigma}\right)-\lambda\text{div}\left(\mathcal{D}(u)\nabla{u}\right)-0.5\mathcal{C}\sigma(u+1)\,,\label{GarLSS163}\\
  \nabla{u}\cdot\mathbf{n}&=\nabla{w}\cdot\mathbf{n}=0\,,\quad \sigma=\sigma_B\quad\textup{on}\ \partial\Omega\,.\label{GarLSS164}
\end{align}
\end{subequations}
The positive constants $\mathcal{P}$, $\mathcal{A}$, and $\mathcal{C}$ denote the proliferation rate, apoptosis rate, and consumption rate, respectively. The mobility is given by $$\mathcal{D}(u)=0.5(1+D)+0.5u(1-D)\,,$$ where $D>0$ is a constant. 

The last two terms on the right-hand side of (\ref{GarLSS161}) represent tumor growth/ proliferation and the process of apoptosis, respectively. The first two terms on the right-hand side of (\ref{GarLSS163}) take chemotaxis and active transport into account. An equivalent term does not appear in the HKNZ model. The last term on the right-hand side of (\ref{GarLSS163}) represents consumption of the nutrient only in the presence of the tumor cell. 

In the spirit of tumor growth modeling, model GLSS is considered with a quasi-steady nutrient (i.e., neglecting the left-hand side of (\ref{GarLSS163})), and Dirichlet boundary conditions are imposed on the nutrient.

\section{Time discretization}
\label{sec:time}
We use the time discretization schemes that have been proposed by the authors of the HKNZ and GLSS models, respectively. In the following, $n\in\mathbb{N}$ denotes the time step and $\tau>0$ the time step size.

\subsection{HKNZ model}
\label{Hilhorst-model-time}
We adopt the recently proposed second-order convexity splitting scheme \cite{WuZZ14} written as
\begin{subequations}\label{scheme:HKNZ}
\begin{gather}
   \frac{u^{(n+1)}-u^{(n)}}{\tau}=\Delta{w^{(n+1)}}+\frac{1}{\varepsilon}\,\tilde p^{(n+\nicefrac{1}{2})}\left(\frac{\sigma^{(n+1)}+\sigma^{(n)}}{2}-s w^{(n+1)}\right)\,,\label{HikKN1_time}
\\
\begin{split}
   w^{(n+1)} &= -\varepsilon\,\Delta{\left(\frac{u^{(n+1)}+u^{(n)}}{2}\right)}+\frac{1}{\varepsilon}\tilde\psi'\left(u^{(n+1)},u^{(n)}\right)\\
             &\qquad\qquad-\alpha_1\tau\Delta{\left(u^{(n+1)}-u^{(n)}\right)}+\alpha_2\tau\left(u^{(n+1)}-u^{(n)}\right)
\end{split}
\label{HikKN2_time}
\\
   \frac{\sigma^{(n+1)}-\sigma^{(n)}}{\tau} = \Delta{\left(\frac{\sigma^{(n+1)}+\sigma^{(n)}}{2}\right)}-\frac{1}{\varepsilon}\,\tilde p^{(n+\nicefrac{1}{2})}\left(\frac{\sigma^{(n+1)}+\sigma^{(n)}}{2}-s w^{(n+1)}\right),\label{HikKN3_time}
\end{gather}
\end{subequations}
with homogeneous Neumann boundary conditions,  where
\begin{displaymath}
   \tilde p^{(n+\nicefrac{1}{2})}=p{\left(\frac{3}{2}u^{(n)}-\frac{1}{2}u^{(n-1)}\right)}\,,
\end{displaymath}
and we take $u^{(-1)}=u^{(0)}$ to initialize the method. Scheme~\eqref{scheme:HKNZ} is a modification of the Crank--Nicolson scheme with a special treatment of the nonlinear terms and two types of stabilizations. In particular, the free energy $\psi(u)$ is split into a convex $\psi_c(u)$ and a concave part $-\psi_e(u)$, which introduces a splitting of $\psi'(u)=\psi'_c(u)-\psi'_e(u)$ that is treated in an implicit and explicit manner, respectively. In detail, we have
\begin{align*}
 \tilde\psi'\left(u^{(n+1)},u^{(n)}\right)  := &\ \psi'_c(u^{(n+1)})-\frac{u^{(n+1)}-u^{(n)}}{2}\psi''_c(u^{(n+1)}) \\
    &-\psi'_e(u^{(n)})-\frac{u^{(n+1)}-u^{(n)}}{2}\psi''_e(u^{(n)})\,.
\end{align*}
\noindent Using the splitting $\psi_c(u)=u^2+\frac{1}{4}$ and $\psi_e(u)=\frac{3}{2}u^2-\frac{1}{4}u^4$ from \cite{WuZZ14}, we 
obtain $$ \tilde\psi'\left(u^{(n+1)},u^{(n)}\right)= u^{(n+1)}\left(\frac{3}{2}\left(u^{(n)}\right)^2-\frac{1}{2}\right)-\frac{1}{2}u^{(n)}-\frac{1}{2}\left(u^{(n)}\right)^3\,,$$ and hence \eqref{scheme:HKNZ} forms a linear system.

The $\alpha_1$-stabilization in \eqref{HikKN2_time} is an artificial diffusivity, whereas the $\alpha_2$-stabilization can be thought of as artificial convexity; see \cite[p.~189]{WuZZ14}. For $\alpha_1,\alpha_2$ sufficiently large, the time-discrete scheme \eqref{scheme:HKNZ} is stable in the sense that the corresponding energy functional decreases with time. We will come back to this system in Section~\ref{Hilhorst-model-space} when we discretize it in space.
 
\subsection{GLSS model}
\label{Garcke-model-time} 
We adopt the time discretization scheme in \cite{GarLSS16}. Using the smooth potential function, we split it into  
$$\psi'(u)=u^3-u\,,$$ 
and treat its nonlinear term implicitly and its linear term explicitly \cite{GarLSS16}. This leads to the following time-discrete scheme in weak form:
\begin{subequations} \label{scheme:GLSS}
\begin{gather}
\begin{split}
  \frac{1}{\tau}\left(u^{(n+1)}-u^{(n)},v\right) =& -\left(\frac{1}{2}\left(1+u^{(n)}\right)^2\nabla{w^{(n+1)}},\nabla{v}\right)  \\
  &+ \left(\left(\mathcal{P}\sigma^{(n)}-\mathcal{A}\right)\left(u^{(n+1)}+1\right),v\right)\,,
\end{split}\label{GarLSS161time1}
\\
\begin{split}
  \left(w^{(n+1)},v\right) = & \beta\varepsilon\left(\nabla{u^{(n+1)}},\nabla{v}\right)+\frac{\beta}{\varepsilon}\left(\left(\left(u^{(n+1)}\right)^3-u^{(n)}\right),v\right) \\
  & -\chi_{u}\left(\sigma^{(n)},v\right)\,,
\end{split}\label{GarLSS162time2}
\\
\begin{split}
  0 =& -\left(\mathcal{D}\left(u^{(n+1)}\right)\nabla{\sigma^{(n+1)}},\nabla{\tilde v}\right)+\lambda\left(\mathcal{D}\left(u^{(n+1)}\right)\nabla{u^{(n+1)}},\nabla{\tilde v}\right)\\
    & -\frac{\mathcal{C}}{2}\left(\sigma^{(n+1)}\left(u^{(n+1)}+1\right),\tilde v\right)\,,
\end{split}\label{GarLSS163time3}
\end{gather}
\end{subequations}
for all $v\in H^{1}(\Omega)$ and $\tilde v\in H_0^{1}(\Omega)$. This scheme decouples the solution for $\sigma^{(n+1)}$, i.e., equation (\ref{GarLSS163time3}), from the remaining system. Hence, we will first solve the nonlinear system (\ref{GarLSS161time1})--(\ref{GarLSS162time2}) to obtain $\left(u^{(n+1)},w^{(n+1)}\right)$. Using this solution in (\ref{GarLSS163time3}) as fixed yields a linear equation that is solved to obtain $\sigma^{(n+1)}$.

Using the nonsmooth obstacle potential instead, (\ref{GarLSS162time2}) will become a variational inequality. To circumvent this difficulty, we proceed as in our previous work \cite{BSB14} and replace the original problem by its Moreau--Yosida regularized version with the energy functional
\begin{displaymath}
\mathcal{E}_{\text{GLSS}}^{\text{reg}}(u,\sigma)=\mathcal{E}_{\text{GLSS}}(u,\sigma)+\frac{1}{2c}\|\max(0,u-1)\|^{2}+\frac{1}{2c}\|\min(0,u+1)\|^{2}.
\end{displaymath}
Here, $0<c\ll 1$ denotes the associated regularization or penalty parameter. As in \cite{GarLSS16}, we treat the nonsmooth potential $\psi_0'(u)=-u$ explicitly, resulting in the regularized scheme
\begin{subequations} \label{scheme:GLSSdisc}
\begin{gather}
\begin{split}
  \frac{1}{\tau}\left(u^{(n+1)}-u^{(n)},v\right) =& -\left(\frac{1}{2}\left(1+u^{(n)}\right)^2\nabla{w^{(n+1)}},\nabla{v}\right) \\
  & + \left(\left(\mathcal{P}\sigma^{(n)}-\mathcal{A}\right)\left(u^{(n+1)}+1\right),v\right)\,,
\end{split}
\label{GarLSS161time1nonsm}
\\
\begin{split}
  \left(w^{(n+1)},v\right) = & \ \beta\varepsilon\left(\nabla{u^{(n+1)}},\nabla{v}\right)-\frac{\beta}{\varepsilon}\left(u^{(n)},v\right)-\chi_{u}\left(\sigma^{(n)},v\right)\\
   &+\frac{1}{c}\left(\max\left(0,u^{(n+1)}-1\right)+\min\left(0,u^{(n+1)}+1\right),v\right)\,,
\end{split}
\label{GarLSS162time2nonsm}
\\
\begin{split}
  0 = &-\left(\mathcal{D}\left(u^{(n+1)}\right)\nabla{\sigma^{(n+1)}},\nabla{\tilde v}\right)+\lambda\left(\mathcal{D}\left(u^{(n+1)}\right)\nabla{u^{(n+1)}},\nabla{\tilde v}\right)\\
    & -\frac{\mathcal{C}}{2}\left(\sigma^{(n+1)}\left(u^{(n+1)}+1\right),\tilde v\right)\,,
\end{split}
\label{GarLSS163time3nonsm}
\end{gather}
\end{subequations}
for all $v\in H^{1}(\Omega)$ and $\tilde v\in H_0^{1}(\Omega)$. As in \eqref{scheme:GLSS}, this scheme decouples the solution for $\sigma^{(n+1)}$, i.e., equation (\ref{GarLSS163time3nonsm}), from the remaining system. In order to ease the notation, from now on we write $u^{\textup{old}}$, $\sigma^{\textup{old}}$, $u$, $w$, and $\sigma$ instead of $u^{(n)}$, $\sigma^{(n)}$, $u^{(n+1)}$, $w^{(n+1)}$, and $\sigma^{(n+1)}$, respectively.

\section{Semi-smooth Newton method}
\label{sec:newton}
We apply the function space-based algorithm motivated in \cite{HinHT11} for solving the regularized time-discrete problem (\ref{GarLSS161time1nonsm})--(\ref{GarLSS162time2nonsm})~\footnote{The nonlinear time-discrete problems using the smooth potential are solved using classical Newton method. We come back to them in the next section when we consider the resulting linear systems after discretizing the problems in space.}. For a specified sequence $c\to 0$, we solve 
the system (\ref{GarLSS161time1nonsm})--(\ref{GarLSS163time3nonsm}), compactly written as
\begin{equation}
  \label{F}
  F_{c}(u_{c}, w_{c})=(F_{c}^{(1)}(u_{c},w_{c}),F_{c}^{(2)}(u_{c},w_{c}))=0\,,
\end{equation}

\noindent for every $c$ by a semi-smooth Newton (SSN) algorithm. 
In (\ref{F}), the components are defined by
\begin{align*}
  \left\langle F_{c}^{(1)}(u,w),v\right\rangle &=
  \frac{1}{\tau}\left(u-u^{\text{old}},v\right) + \left(\frac{1}{2}\left(1+u^{\text{old}}\right)^2\nabla{w},\nabla{v}\right)\\
    &\quad -\left(\left(\mathcal{P}\sigma^{\text{old}}-\mathcal{A}\right)\left(u+1\right),v\right) ,\\
  \left\langle F_{c}^{(2)}(u,w),v\right\rangle &=
  \left(w,v\right) - \beta\varepsilon\left(\nabla{u},\nabla{v}\right) + \frac{\beta}{\varepsilon}\left(u^{\text{old}},v\right) + \chi_{u}\left(\sigma^{\text{old}},v\right)\\
     &\quad -\frac{1}{c}\left(\max\left(0,u-1\right)+\min\left(0,u+1\right),v\right) ,
  \end{align*}
\noindent for all $u,w,v\in H^{1}(\Omega)$. 
 
Due to the presence of the max- and min-operators, $F_{c}$ is not Fr\'echet-differentiable. However, it satisfies the weaker notion of Newton differentiability; see \cite[p.~866--867]{HinIK03} for the definition as well as a superlinear convergence result for the (semi-smooth) Newton iteration
\begin{equation*}
  \label{semi-newton}
  x^{(k+1)}=x^{(k)}-G(x^{(k)})^{-1}F(x^{(k)}),\quad k=0,1,\ldots
\end{equation*}
with $G$ being a Newton-derivative of $F$.

In analogy to \cite[p.~788]{HinHT11} and \cite[pp.~885--886]{HinIK03}, the operator $G_{c}(u,w)$ given by
\begin{multline*}
  \left\langle G_{c}(u,w)(\delta{u},\delta{w}),(\phi,\psi)\right\rangle=\\
  \left(\begin{array}{c}
     \frac{1}{\tau}\left(\delta{u},\phi\right) + \left(\frac{1}{2}\left(1+u^{\text{old}}\right)^2\nabla{\delta{w}},\nabla{\phi}\right)
          -\left(\left(\mathcal{P}\sigma^{\text{old}}-\mathcal{A}\right)\delta{u},\phi\right)\\
     \left(\delta{w},\psi\right) - \beta\varepsilon\left(\nabla{\delta{u}},\nabla{\psi}\right)
          -\frac{1}{c}(\chi_{\mathcal{M}(u)}\delta{u},\psi)
  \end{array}\right)
\end{multline*}
serves as a Newton-derivative for $F_{c}$, where $\chi_{\mathcal{M}(u)}$ is the characteristic function of the set $\mathcal{M}(u)\coloneqq\{x\in\Omega : |u({\bf x})|> 1\}\,.$

\section{Finite element approximation}\label{sec:FEM}
In this section, we apply the finite element method \cite{StrF73} to the regularized system in \eqref{scheme:GLSSdisc}. We will also apply it to the systems that use a smooth potential. Since both procedures are similar, we only present the methodology based on the regularized setting. For the smooth cases, we will state the fully discrete linear system at the end of this section.

\subsection{Discretization of the regularized GLSS model (\ref{GarLSS161time1nonsm})--(\ref{GarLSS162time2nonsm})}
In the following, we assume for simplicity that $\Omega$ is a polyhedral domain. Generalizations to curved domains are possible using boundary finite elements with curved faces. Let $\left\{\mathcal{R}_{h}\right\}_{h>0}$ be a triangulation of $\Omega$ with maximal element size $h$.  Let $J_{h}$ be the set of nodes of $\mathcal{R}_{h}$ and let $p_{j}\in J_{h}$ be the coordinates of these nodes. We approximate the infinite-dimensional space $H^{1}(\Omega)$ by the finite-dimensional space
\begin{displaymath}
 S_{h}\coloneqq\{\phi\in C^{0}(\overline{\mskip-.5\thinmuskip \Omega\mskip-.7\thinmuskip}):\phi\left.\right|_{R}\in Q_{1}(R)\ \ \forall R\in\mathcal{R}_{h}\}\subset H^{1}(\Omega)
\end{displaymath}
of continuous, piecewise multilinear functions. For instance, in the two-dimensional case $d=2,$ we use bilinear functions, i.e., ${Q}_{1}=\textup{span}\{x^{\alpha_{i}}y^{\alpha_{i}}: \alpha_{i}\in\{0,1\},\, i=1,2\}$. We denote the standard nodal basis functions of $S_{h}$ by $\varphi_{j}$ for all $j\in J_{h}$. They have the property $\varphi_{j}(p_{i})=\delta_{ij},\, i,j=1,\ldots,m$. The discretized version of the problem (\ref{GarLSS161time1nonsm})--(\ref{GarLSS162time2nonsm}) is the following: 

\noindent Given $u_{h}^{\textup{old}}\in S_{h}$ and 
$\sigma_{h}^{\textup{old}}\in S_{h}^B:=\left\{s_h\in S_h \left|\right. s_h=\sigma_B\ \text{on}\ \partial\Omega\right\}$, 
find $(u_{c,h},w_{c,h})\in S_{h}\times S_{h}$ such that
\begin{subequations}\label{newton-eq-diskret}
\begin{gather}
\left\langle F_{c,h}^{(1)}(u_{c,h},w_{c,h}),v_{h}\right\rangle=0\quad\forall v_{h}\in S_{h}\,, \label{newton-eq-diskret1} \\
\left\langle F_{c,h}^{(2)}(u_{c,h},w_{c,h}),v_{h}\right\rangle=0\quad\forall v_{h}\in S_{h}\,,\label{newton-eq-diskret2}
\end{gather}
\end{subequations}
\noindent and the components are
\begin{align*}
  \left\langle F_{c,h}^{(1)}(u_{c,h},w_{c,h}),v_h\right\rangle =&
  \frac{1}{\tau}\left(u_{c,h}-u_h^{\text{old}},v_h\right)_h 
    + \left(\frac{1}{2}\left(1+u_h^{\text{old}}\right)^2\nabla{w_{c,h}},\nabla{v_h}\right)\\
    & -\left(\left(\mathcal{P}\sigma_h^{\text{old}}-\mathcal{A}\right)\left(u_{c,h}+1\right),v_h\right)_h\,,\\
  \left\langle F_{c,h}^{(2)}(u_{c,h},w_{c,h}),v_h\right\rangle =&
  \left(w_{c,h},v_h\right)_h - \beta\varepsilon\left(\nabla{u_{c,h}},\nabla{v_h}\right) 
     + \frac{\beta}{\varepsilon}\left(u_h^{\text{old}},v_h\right)_h 
     \\
     & + \chi_{u}\left(\sigma_h^{\text{old}},v_h\right)_h\\
     & -\frac{1}{c}\left(\max\left(0,u_{c,h}-1\right)+\min\left(0,u_{c,h}+1\right),v_h\right)_h\,.
  \end{align*}

The semi-inner product $(\cdot,\cdot)_{h}$ on $C_{0}(\overline{\mskip-.5\thinmuskip \Omega\mskip-.7\thinmuskip})$ is defined by
\begin{displaymath}
	(f,g)_{h}\coloneqq\int_{\Omega}{\pi_{h}(f({\bf x})g({\bf x}))\D{{\bf x}}}=\sum_{i=1}^{m}{(1,\varphi_{i})f(p_{i})g(p_{i})}\quad\forall f,g\in C_{0}(\overline{\mskip-.5\thinmuskip \Omega\mskip-.7\thinmuskip})\,,
\end{displaymath}

\noindent where $\pi_{h}\colon C_{0}(\overline{\mskip-.5\thinmuskip \Omega\mskip-.7\thinmuskip})\to S_{h}$ is the Lagrange interpolation operator. Within our finite element framework, for a given every step of the SSN method for solving \eqref{newton-eq-diskret} requires to compute $(\delta{u_{h}},\delta{w}_{h})\in S_{h}\times S_{h}$ satisfying
\begin{subequations}
\begin{gather}
\begin{split}
      \frac{1}{\tau}\left(\delta{u_h},v_h\right)_h + \left(\frac{1}{2}\left(1+u_h^{\text{old}}\right)^2\nabla{\delta{w_h}},\nabla{v_h}\right)
          &-\left(\left(\mathcal{P}\sigma_h^{\text{old}}-\mathcal{A}\right)\delta{u_h},v_h\right)_h \\
            &= -F_{c,h}^{(1)}(u_{h},w_{h})\,,
\end{split}
\label{eqfem1}
\\
\begin{split}
     \left(\delta{w_h},v_h\right)_h - \beta\varepsilon\left(\nabla{\delta{u_h}},\nabla{v_h}\right)
          &-\frac{1}{c}(\chi^h_{\mathcal{M}(u_h)}\delta{u_h},v_h)_h \\
            &= -F_{c,h}^{(2)}(u_{h},w_{h})\,,
\end{split}
\label{eqfem2}
\end{gather}
\end{subequations}
\noindent for all $v_{h}\in S_{h}$, where $\chi_{\mathcal{M}(u_{h})}^{h}\coloneqq\sum_{i=1}^{m}{\chi_{\mathcal{M}(u_{h})}^{h}(p_{i})\varphi_{i}}$ with $\chi_{\mathcal{M}(u_{h})}^{h}(p_{i})=0$ if $-1\leq u_{h}(p_{i})\leq 1$ and $\chi_{\mathcal{M}(u_{h})}^{h}(p_{i})=1$ otherwise. If we now write a function $v_{h}\in S_{h}$ by $v_{h}=\sum_{j\in J_{h}}{v_{h,j}\,\varphi_{j}}$ and denote the vector of coefficients by ${\bf v}$, the fully decoupled discrete linear system per SSN step reads in matrix form as 
\begin{equation}
\begin{split}
\label{GLSS_nonsmooth_matrix_1}
\left[\begin{array}{cc}
 \bm{M} & \frac{1}{c} \bm{G} + \beta\varepsilon \bm{K}  \\
  -\frac{\tau}{2} \bm{\hat K} & \left[\bm{I}_m-\tau\left(\mathcal{P}\bm{\sigma}^{\text{old}}-\mathcal{A}\bm{I}_m\right)\right]\bm{M}
\end{array}\right] 
\left[\begin{array}{c} -\bm{w}^{(k+1)}\\ \bm{u}^{(k+1)}\end{array}\right]\\
 = \left[\begin{array}{c}
  \frac{1}{c}\left(\bm{G}_+ - \bm{G}_-\right)\bm{1}
    + \frac{\beta}{\varepsilon}\bm{M}\bm{u}^{\text{old}}
    + \chi_u\bm{M}\bm{\sigma}^{\text{old}}\\
  \bm{M} \bm{u}^{\textup{old}} + \tau\left(\mathcal{P}\bm{\sigma}^{\text{old}}-\mathcal{A}\right)\bm{M}\bm{1}
 \end{array}\right]\,.
\end{split}
\end{equation}
Here, $\bm{u}^{(k+1)}, \bm{w}^{(k+1)}\in\mathbb{R}^{m}$ are the unknowns of the current SSN step and $\bm{u}^{\textup{old}},\bm{\sigma}^{\textup{old}}\in\mathbb{R}^{m}$ are the solution vectors from the previous time step. 

The presence of $\bm{\sigma}^{\textup{old}}$ in the coefficient matrix has the meaning of an $m\times m$ diagonal matrix containing the entries of $\bm{\sigma}^{\textup{old}}$ on its diagonal. $\bm{I}_m\in\mathbb{R}^{m\times m}$ and $\bm{1}$ denotes the identity matrix and the vector of ones, respectively. The lumped mass matrix $\bm{M}$, stiffness matrix $\bm{K}$, and modified stiffness matrix $\bm{\hat K}$ are defined as
\begin{align*}
\bm{M}&\coloneqq((\varphi_{j},\varphi_{i})_{h})_{i,j=1,\ldots,m}=\textup{diag}((1,\varphi_{i}))_{i=1,\ldots,m}\in\mathbb{R}^{m\times m}\,,\\
\bm{K}&\coloneqq((\nabla{\varphi_{j}},\nabla{\varphi_{i}}))_{i,j=1,\ldots,m}\in\mathbb{R}^{m\times m}\,,\\
\bm{\hat K}&\coloneqq\left(\left(\left[1+u_h^{\text{old}}\right]^2\nabla{\varphi_{j}},\nabla{\varphi_{i}}\right)\right)_{i,j=1,\ldots,m}\in\mathbb{R}^{m\times m}\,.
\end{align*}

The matrix representations coming from the generalized derivative of the term $(\chi_{\mathcal{M}(u_{h})}^{h}\delta{u}_{h},v_{h})_{h}$ are the following diagonal matrices
\begin{align*}
 \bm{G}&=\bm{G}(\bm{u}^{(k)})=\textup{diag}\left(\begin{array}{rl} [\bm{M}]_{ii}&\textup{if}\ |u^{(k)}_{h,i}|>1,\\0&\textup{otherwise}\end{array}\right)_{i=1,\ldots,m}\in\mathbb{R}^{m\times m},\\
 \bm{G}_{+}&=\bm{G}_{+}(\bm{u}^{(k)})=\textup{diag}\left(\begin{array}{rl}[\bm{M}]_{ii}&\textup{if}\ u^{(k)}_{h,i}>1,\\0,&\textup{otherwise}\end{array}\right)_{i=1,\ldots,m}\in\mathbb{R}^{m\times m},\\
 \bm{G}_{-}&=\bm{G}_{-}(\bm{u}^{(k)})=\textup{diag}\left(\begin{array}{rl}[\bm{M}]_{ii}&\textup{if}\ u^{(k)}_{h,i}<-1\\0&\textup{otherwise}\end{array}\right)_{i=1,\ldots,m}\in\mathbb{R}^{m\times m},
\end{align*}
where $\bm{u}^{(k)}\in\mathbb{R}^{m}$ is the solution from the previous SSN step $k$.

\subsection{Discretization of the remaining part (\ref{GarLSS163time3nonsm}) of the GLSS model}
\label{equ3-discrete}
Once we have the solution of \eqref{newton-eq-diskret}, i.e., an approximation to $(u_h,w_h)$, we solve the discretized version of (\ref{GarLSS163time3nonsm}): Given $u_{h}\in S_{h}$, find $\sigma_{h}\in S_{h}^B$ such that
\begin{eqnarray*}
    \frac{1}{2}\left(\left[1+D+u_{h}(1-D)\right]\nabla{\sigma_h},\nabla{\tilde v_h}\right)
     + \frac{\mathcal{C}}{2}\left(\left(u_{h}+1\right)\sigma_{h},\tilde v_h\right)_h \\
      = \frac{\lambda}{2}\left(\left[1+D+u_{h}(1-D)\right]\nabla{u_{h}},\nabla{\tilde v_h}\right)
\end{eqnarray*}
for all $\tilde v_{h}\in S_{h}^0:=\left\{s_h\in S_h \left|\right. s_h=0\ \text{on}\ \partial\Omega\right\}$. Let us denote the standard nodal basis functions of $S_{h}^0$ by $\{\tilde \varphi_{1},\ldots,\tilde \varphi_{m}\}$. In order to satisfy the Dirichlet conditions in $S_{h}^B$, we proceed as in \cite[pp.~17--18]{ElmSW05}: We extend the basis set by defining additional functions $\tilde \varphi_{m+1},\ldots,\tilde \varphi_{m+m_{\partial}}$ and select fixed coefficients $\sigma_{h,j}$, $j=m+1,\ldots,m+m_{\partial}$, so that the function $\sum_{j=m+1}^{m+m_{\partial}}{\sigma_{h,j}\,\tilde \varphi_{j}}$ interpolates the boundary data $\sigma_B$ on $\partial\Omega$. The finite element approximation $\sigma_{h}\in S_{h}^B$ is then uniquely associated with the vector 
$\bm{\sigma}=\left[\begin{array}{ccc}\sigma_{h,1},\ldots,\sigma_{h,m}\end{array}\right]$ of real coefficients in the expansion
$$\sigma_{h}=\sum_{j=1}^{m}{\sigma_{h,j}\,\tilde \varphi_{j}}+\sum_{j=m+1}^{m+m_{\partial}}{\sigma_{h,j}\,\tilde \varphi_{j}}\,.$$

The fully discrete linear system reads in matrix form as 
\begin{equation}
 \label{GLSS_nonsmooth_matrix_2}
 \left[\bm{\hat K}_0+\frac{\mathcal{C}}{2}\left(\bm{I}_m+\bm{u}\right)\bm{M}_0\right]\bm{\sigma}=\bm{g}\,.
\end{equation}
Similarly, the presence of $\bm{u}$ -- which is the coefficient vector of the solution $u_h$ of (\ref{newton-eq-diskret1})--(\ref{newton-eq-diskret2}) -- in the coefficient matrix has the meaning of an $m\times m$ diagonal matrix containing the entries of $\bm{u}$ on its diagonal. The lumped mass matrix $\bm{M}_0$, modified stiffness matrix $\bm{\hat K}_0$, and right-hand side $\bm{g}$ are defined as 
\begin{align*}
\bm{M}_0&\coloneqq((\tilde\varphi_{j},\tilde\varphi_{i})_{h})_{i,j=1,\ldots,m}=\textup{diag}((1,\tilde\varphi_{i}))_{i=1,\ldots,m}\in\mathbb{R}^{m\times m}\,,\\
\bm{\hat K}_0&\coloneqq\left(\left(\frac{1}{2}\left[1+D+u_h\left(1-D\right)\right]\nabla{\tilde\varphi_{j}},\nabla{\tilde\varphi_{i}}\right)\right)_{i,j=1,\ldots,m}\in\mathbb{R}^{m\times m}\,,\\
\bm{g}&\coloneqq\left(\left(\frac{\lambda}{2}\left[1+D+u_h\left(1-D\right)\right]\nabla{u_h},\nabla{\tilde\varphi_{i}}\right)\right)_{i=1,\ldots,m}\\
&\quad -\sum_{j=m+1}^{m+m_{\partial}}{\sigma_{h,j}\left(\left(\frac{1}{2}\left[1+D+u_h\left(1-D\right)\right]\nabla{\tilde\varphi_{j}},\nabla{\tilde\varphi_{i}}\right)\right)_{i=1,\ldots,m}}\in\mathbb{R}^{m}\,.
\end{align*}

Note that $\bm{M}$ and $\bm{M}_0$ are diagonal, symmetric positive definite matrices. Furthermore,  $\bm{\hat K}_0$ is symmetric positive definite, and $\bm{K}$ and $\bm{\hat K}$ are symmetric positive semidefinite. In particular, they have the following eigenvalue characterization.
\begin{prop}[{\cite[pp.~57--60]{ElmSW05}}]
\label{chap3:prop:mass_stiffness}
Let $d\in\{2,3\}$ be the spatial dimension. Then,
\begin{gather*}
	\tilde ch^{d}\leq\frac{(\bm{M}\bm{v},\bm{v})}{(\bm{v},\bm{v})}\leq Ch^{d},\\
	\tilde c_0h^{d}\leq\frac{(\bm{M}_0\bm{v},\bm{v})}{(\bm{v},\bm{v})}\leq C_0h^{d},\\
	0\leq\frac{(\bm{K}\bm{v},\bm{v})}{(\bm{v},\bm{v})}\leq \tilde Ch^{d-2},
\end{gather*}
for all $\bm{v}\in\mathbb{R}^{m}$. The constants $\tilde c,\tilde c_0,C,C_0,\tilde C$ are positive and independent of $h$. In particular, $\bm{K}$ has a one-dimensional kernel spanned by the constant vector $\bm{1}=[1,\ldots,1]^{T}\in\mathbb{R}^{m}$.
\end{prop}

In terms of the condition number, $\kappa(\bm{M})\leq \tilde c^{-1}C$, $\kappa(\bm{M}_0)\leq \tilde c_0^{-1}C_0$, and $\kappa(\bm{K})=\infty$. Note that the original version of Proposition \ref{chap3:prop:mass_stiffness} in \cite[pp.~57--60]{ElmSW05} is actually stated under proper mesh subdivisions. We will not go into these details but note that they hold true in our setting.

Since the system~\eqref{GLSS_nonsmooth_matrix_2} is only a second order problem and is decoupled from the remaining variables, we employ an algebraic multigrid \texttt{AMG} approximation of the coefficient matrix $\bm{\hat K}_0+\frac{\mathcal{C}}{2}\left(\bm{I}_m+\bm{u}\right)\bm{M}_0$. More precisely, a smoothed aggregation multigrid preconditioner~\cite{pyamg2023} constructed and applied with a \texttt{MINRES} solver~\cite{minres}.

\subsection{Discretization of the smooth GLSS model (\ref{GarLSS161time1})--(\ref{GarLSS162time2})}
Using the smooth potential, the fully discrete linear system per Newton step $k$ of the GLSS model reads in matrix form as 
\begin{align}
\begin{split}
\left[\begin{array}{cc}
  \frac{3\beta}{\varepsilon}\left(\bm{u}^{(k)}\right)^2 \bm{M} + \beta\varepsilon \bm{K}  & \bm{M} \\
  \left[\bm{I}_m-\tau\left(\mathcal{P}\bm{\sigma}^{\text{old}}-\mathcal{A}\bm{I}_m\right)\right]\bm{M} & -\frac{\tau}{2} \bm{\hat K}
\end{array}\right]
\left[\begin{array}{c} \bm{u}^{(k+1)} \\ -\bm{w}^{(k+1)}\end{array}\right]\\
=\left[\begin{array}{c}
  \frac{2\beta}{\varepsilon}\bm{M}\left(\bm{u}^{(k)}\right)^3
    + \frac{\beta}{\varepsilon}\bm{M}\bm{u}^{\text{old}}
    + \chi_u\bm{M}\bm{\sigma}^{\text{old}}\\
  \bm{M} \bm{u}^{\textup{old}} + \tau\left(\mathcal{P}\bm{\sigma}^{\text{old}}-\mathcal{A}\bm{I}_m\right)\bm{M}\bm{1}
 \end{array}\right]\,.
\end{split} \label{GLSS_smooth_matrix}
\end{align}
Once we have an approximation to $(u_h,w_h)$, we solve the discretized version of (\ref{GarLSS163time3}) as in Section \ref{equ3-discrete}.

\subsection{Discretization of the smooth HKNZ model (\ref{HikKN1_time})--(\ref{HikKN3_time})} \label{Hilhorst-model-space}
Using the smooth potential, the fully discrete linear system of the HKNZ model reads in matrix form as 
\begingroup
\setlength{\arraycolsep}{2pt}
\begin{equation}
\label{matrix-HKNZ-smooth}
\begin{split}
&\scalebox{0.9}{$\displaystyle
\left[\begin{array}{ccc|c}
\frac{s^2}{\varepsilon}\bm{P}_u\bm{M}+s\bm{K} & &
-\frac{s}{2\varepsilon}\bm{P}_u\bm{M} &
\frac{s}{\tau}\bm{M} \BBstrut\\
-\frac{s}{2\varepsilon}\bm{P}_u\bm{M} & &
\frac{1}{2}\bm{M}\left(\frac{1}{\tau}\bm{I}_m
+\frac{1}{2\varepsilon}\bm{P}_u\right)+\frac{1}{4}\bm{K} &
\bm{0}\\ \hline
\frac{s}{\tau}\bm{M} & & \bm{0} &
-s\bm{M}\left(\frac{1}{\varepsilon\tau}\bm{D}_u+\alpha_2\bm{I}_m\right)
-s\bm{K}\left(\frac{\varepsilon}{2\tau}+\alpha_1\right) \Tstrut
\end{array}\right]
$}
\\
&\qquad \quad\scalebox{0.9}{$\displaystyle
\left[\begin{array}{c}
\bm{w}^{(n+1)}\\
\bm{\sigma}^{(n+1)}\\
\bm{u}^{(n+1)}
\end{array}\right]
=
\left[\begin{array}{c}
\frac{s}{\tau}\bm{M}\bm{u}^{(n)}
+\frac{s}{2\varepsilon}\bm{P}_u\bm{M}\bm{\sigma}^{(n)}
\\
\frac{1}{2}\left(\frac{1}{\tau}\bm{I}_m
-\frac{1}{2\varepsilon}\bm{P}_u\right)
\bm{M}\bm{\sigma}^{(n)}
-\frac{1}{4}\bm{K}\bm{\sigma}^{(n)}
\\
s\left(\frac{\varepsilon}{2\tau}-\alpha_1\right)\bm{K}\bm{u}^{(n)}
-s\left(\frac{1}{2\varepsilon\tau}+\alpha_2\right)\bm{M}\bm{u}^{(n)}
-\frac{s}{2\varepsilon\tau}\bm{M}\left(\bm{u}^{(n)}\right)^3
\end{array}\right],
$}
\end{split}
\end{equation}
\endgroup
where $\bm{P}_u\in\mathbb{R}^{m\times m}$ is a diagonal matrix representing $p{\left(\frac{3}{2}u^{(n)}-\frac{1}{2}u^{(n-1)}\right)}$ pointwise. $\bm{D}_u\in\mathbb{R}^{m\times m}$ is a diagonal matrix representing 
$\frac{3}{2}\left(u^{(n)}\right)^2-\frac{1}{2}$, i.e., the implicit part of $\tilde\psi'\left(u^{(n+1)},u^{(n)}\right)$, pointwise.

\section{Preconditioning}\label{sec:prec}
In this section, we focus on the development of the linear systems arising from the discretized tumor models.
\subsection{Preconditioning for the smooth HKNZ model (\ref{scheme:HKNZ})}\label{sec:prec_smoothH}
Looking at the estimated magnitude of entries of the blocks in the coefficient matrix of (\ref{matrix-HKNZ-smooth}), we have that 
$\mathcal{O}\left(\frac{s}{2\varepsilon}\bm{P}_u\bm{M}\right)=s\varepsilon^{d-1}$ since $h=\mathcal{O}(\varepsilon)$. 
Hence, we propose to neglect this block within our approximation. Further, we typically have $\tau=\mathcal{O}(\varepsilon^2)$ in numerical simulations, which leads us to $\mathcal{O}\left(\frac{1}{2}\bm{M}\left(\frac{1}{\tau}\bm{I}_m+\frac{1}{2\varepsilon}\bm{P}_u\right)+\frac{1}{4}\bm{K}\right) =\frac{1}{2\tau}\bm{M}+\frac{1}{4}\bm{K}$. All in all, we suggest the following preconditioners
\begin{align*}
\bm{\mathcal{P}}_{diag} &=\left[\begin{array}{cccc|c}
\frac{s^2}{\varepsilon}\bm{M}+s\bm{K} & & & \bm{0} & \bm{0} \BBstrut\\
\bm{0} & & & \frac{1}{2\tau}\bm{M}+\frac{1}{4}\bm{K} & \bm{0}\\\hline
\bm{0} & & & \bm{0} & -\Shat \Tstrut
\end{array}\right]=:
\left[\begin{array}{ccc|c}
\bm{A}_1 & &  \bm{0} & \bm{0} \Bstrut\\
\bm{0} & &  \bm{A}_2 & \bm{0}\\\hline
\bm{0} & & \bm{0} & -\Shat \Tstrut
\end{array}\right]\,, \\
\bm{\mathcal{P}}_{tri}&=\left[\begin{array}{cccc|c}
\frac{s^2}{\varepsilon}\bm{M}+s\bm{K} & & & \bm{0} & \bm{0} \BBstrut\\
\bm{0} & & & \frac{1}{2\tau}\bm{M}+\frac{1}{4}\bm{K} & \bm{0}\\\hline
\frac{s}{\tau}\bm{M} & & & \bm{0} & -\Shat \Tstrut
\end{array}\right]=:
\left[\begin{array}{ccc|c}
\bm{A}_1 & &  \bm{0} & \bm{0} \Bstrut\\
\bm{0} & &  \bm{A}_2 & \bm{0}\\\hline
\frac{s}{\tau}\bm{M} & & \bm{0} & -\Shat \Tstrut
\end{array}\right]\,.
\end{align*}

The formula for the exact Schur complement is rather involved and we omit its explicit presentation. We drop one of its terms to obtain
\begin{displaymath}
\bm{S}\approx
 -s\bm{M}\left(\frac{1}{\varepsilon\tau}\bm{D}_u+\alpha_2\bm{I}_m\right)-s\bm{K}\left(\frac{\varepsilon}{2\tau}+\alpha_1\right) -
\frac{s^2}{\tau^2}\bm{M}\bm{A}_1^{-1}\bm{M},
\end{displaymath}
and further approximate it by the practical Schur complement approximation
\begin{align*}
 \Shat &= -\left(\sqrt{\frac{\varepsilon}{2\tau}}\bm{A}_1+\frac{s}{\tau}\bm{M}\right)\bm{A}_1^{-1}\left(\sqrt{\frac{\varepsilon}{2\tau}}\bm{A}_1+\frac{s}{\tau}\bm{M}\right).
\end{align*}

In practice, the application of $\bm{\mathcal{P}}^{-1}$ requires the approximate inversion of three operators: the blocks $\bm{A_1}$, $\bm{A_2}$ and the factor of $\hat{\bm{S}}$, namely $\left(\sqrt{\frac{\varepsilon}{2\tau}}\bm{A}_1+\frac{s}{\tau}\bm{M}\right)$. For the approximate inversion of the diagonal blocks, we employ algebraic multigrid (\texttt{AMG}) using the Ruge--St\"{u}ben strategy~\cite{pyamg2023}. Since the HKNZ formulation leads to a linear system at each time step, no outer nonlinear iteration is required and the system is solved directly.

\subsection{Preconditioning for the smooth GLSS model (\ref{GarLSS161time1})--(\ref{GarLSS162time2})}\label{sec:prec_smoothG}
We suggest the lower triangular block preconditioner 
\begin{equation}
\bm{\mathcal{P}}=\left[\begin{array}{cc}
  \frac{3\beta}{\varepsilon}\left(\bm{u}^{(k)}\right)^2 \bm{M} + \beta\varepsilon \bm{K}  & \bm{0} \\
  \left[\bm{I}_m-\tau\left(\mathcal{P}\bm{\sigma}^{\text{old}}-\mathcal{A}\bm{I}_m\right)\right]\bm{M} & -\bm{S}
\end{array}\right] =: \left[\begin{array}{cc}
   \bm{A} & \bm{0} \\
  \bm{B} & -\bm{S}
\end{array}\right]
\end{equation}
where $\bm{S} = - \left(\frac{\tau}{2} \bm{\hat K} +  \bm{B} \bm{A}^{-1} \bm{M}\right) $ is the Schur complement. We approximate the Schur complement
\begin{displaymath}
\bm{S}\approx \bm{\hat S} = -\left(\frac{\sqrt{\tau}}{2} \bm{\hat K} + \bm{B} \right)  \bm{A}^{-1} \left( \sqrt{\tau} \bm{A} +\bm{M}\right)
\end{displaymath}

Here, the inversion of three operators, the block $\bm{A}$, and the two factors of $\hat{\bm{S}}$, namely $\frac{\sqrt{\tau}}{2}\hat{\bm{K}} + \bm{B}$ and $\sqrt{\tau}\bm{A} + \bm{M}$ are approximated by employing incomplete LU (\texttt{ILU}) factorizations with a fill factor of $10$ and a drop tolerance of $10^{-4}$ in two dimensions, and  the \texttt{AMG} in three dimensions; see Section~\ref{sec:results} for a comparison of both approaches. 

\subsection{Preconditioning for the nonsmooth GLSS model (\ref{GLSS_nonsmooth_matrix_1})--(\ref{GLSS_nonsmooth_matrix_2})}\label{sec:prec_ns}
Here, we again suggest the lower triangular block preconditioner
\begin{equation}
\bm{\mathcal{P}}= \left[\begin{array}{cc}
 \bm{M}   & \bm{0} \\
-\frac{\tau}{2} \bm{\hat K} & \bm{S}
\end{array}\right] 
\end{equation}
where $\bm{S} =  \left[\bm{I}_m- \tau\left(\mathcal{P}\bm{\sigma}^{\text{old}}-\mathcal{A}\bm{I}_m\right)\right]\bm{M} +  \frac{\tau}{2} \bm{\hat K} \bm{M}^{-1}\left(\frac{1}{c} \bm{G} + \beta\varepsilon \bm{K}\right) =  \bm{\tilde M} + \frac{\tau}{2} \bm{\hat K} \bm{M}^{-1} \bm{\tilde K}$ is the Schur complement. We approximate the Schur complement
\begin{displaymath}
\bm{S}\approx \bm{\hat S} = \left(\frac{\sqrt{\tau}}{2} \bm{\hat K} + \bm{\tilde M} \right)  \bm{M}^{-1} \left( \sqrt{\tau} \bm{\tilde K} +\bm{M}\right).
\end{displaymath}
As in the smooth case, the diagonal blocks are inverted approximately using \texttt{ILU} factorizations.

\section{Numerical results} \label{sec:results}
In this section, we present numerical results for the Cahn--Hilliard tumor growth model, considering two formulations:  the smooth potential of HKNZ model in Section~\eqref{sec:Hilhorst-model} and the (non)smooth potential and of GLSS model in Section~\eqref{sec:Garcke-model}. For both formulations, the block preconditioners we introduce employ left preconditioning and can be incorporated into a range of Krylov subspace solvers.

In HKNZ model, the symmetric saddle-point structure of this formulation admits the use of \texttt{MINRES}~\cite{minres} as the Krylov subspace solver, with a relative tolerance of $\texttt{rtol} = 10^{-10}$ for the preconditioned relative residual. 

In GLSS model, for the nonsymmetric saddle-point system arising from the Newton linearization, we employ \texttt{Bi-CGSTAB}~\cite{bicgstab}, although other solvers such as \texttt{GMRES}~\cite{gmres} or \texttt{QMR}~\cite{QMR} are equally compatible with the preconditioner presented here. The \texttt{Bi-CGSTAB} relative tolerance is set to $\texttt{rtol} = 10^{-15}$ for the preconditioned relative residual, following the \texttt{scipy.sparse.linalg} convention~\cite{scipy}. For the (smooth and nonsmooth) Newton method, we use the stopping criterion in~\cite{HinHT11}, given by
\begin{equation}
  \bigl\|F(\nu_{h}^{(k)})\bigr\|_{2}
  \leq
  \epsilon_{\mathrm{rel}}\,\bigl\|F(\nu_{h}^{(0)})\bigr\|_{2}
 + \epsilon_{\mathrm{abs}},
  \qquad k = 1, \ldots, k_{\mathrm{max}},
\end{equation}
where $F(\nu_h^{(k)})$ denotes the residual of the nonlinear system evaluated at the $k$-th Newton iterate $\nu_h^{(k)} = (u_h^{(k)}, w_h^{(k)})$, and $\epsilon_{\mathrm{rel}}$, $\epsilon_{\mathrm{abs}} > 0$ are prescribed relative and absolute tolerances, respectively. We set $k_{\mathrm{max}} = 20$, $\epsilon_{\mathrm{rel}} = 10^{-12}$ and $\epsilon_{\mathrm{abs}} = 10^{-6}$ in all experiments. Unlike the obstacle formulation, the smooth double-well potential requires no continuation in a regularization parameter; the Newton method is applied directly with the full nonlinearity at every time step, using the solution from the previous time step as the initial guess.

Both formulations are discretized using linear finite elements on triangulated meshes. All implementations are carried out using the finite element package \texttt{DOLFINx}~\cite{dolfinx}, the next-generation interface to the \texttt{FEniCS} project~\cite{fenics}.

For spatial adaptivity, applied to the Garcke formulation only, we follow the mesh refinement strategy of~\cite{barrett2004finite, garcke2022viscoelastic, GarLSS16, garcke2022numerical}: at each time step we mark all elements $K \in \mathcal{T}^{h}$ satisfying $|u_{h}^{n-1}(x)| \leq 1 - \xi$ for some $\xi > 0$, together with their immediate neighbors, provided $h_{K} = \operatorname{diam}(K) > h_{f}$. Since \texttt{DOLFINx} does not currently support local mesh coarsening, we reset the mesh to a uniform coarse triangulation of size $h_{c}$ at the start of each refinement step and reapply the local refinement procedure. The solution is transferred directly from the previous fine mesh to the newly refined mesh. The HKNZ formulation is solved on a fixed uniform mesh throughout, as the smoother potential does not require the fine interface resolution provided by adaptivity.

\subsection{Model HKNZ}
In this section, we will provide some numerical results for the linear system~\eqref{matrix-HKNZ-smooth}.  
\begin{table}[htp!]
\centering
\caption{Parameter values used for the HKNZ model experiments.}
\label{tab:hilhorts_params}
\begin{tabular}{c|cccc}
\hline
In 2D & $\Omega=[-1,1]^2$  & $n_{x_1} \times n_{x_2} = 128 \times 128$ & $\varepsilon= 0.015$  \\
In 3D & $\Omega=[-1,1]^3 $ & $n_{x_1} \times n_{x_2} \times n_{x_3} = 64 \times 64 \times 64 $   & $\varepsilon= 0.03$  \\
\hline
   & $\tau=1\cdot 10^{-5}$& $p_0= 3$ & $s = 0.01$ \\
   & $\alpha_1=\alpha_2=2$ & $\sigma_0(x) = 1$ & \\
\hline
\end{tabular}
\end{table}

\subsubsection{Single tumor growth simulation}
\label{subsec:singel_tumor}
The initial condition is a single spherical tumor centered at the origin,
\begin{equation}
  u_0(x_1, x_2)
  = \tanh\!\left(\frac{r_0 - r}{\sqrt{2}\,\varepsilon}\right),
  \qquad
  r = \sqrt{x_1^2 + x_2^2},
\end{equation}
with radius $r_0 = 0.15$.

\begin{figure}[htp!]
    \centering
    \includegraphics[width=0.88\linewidth]{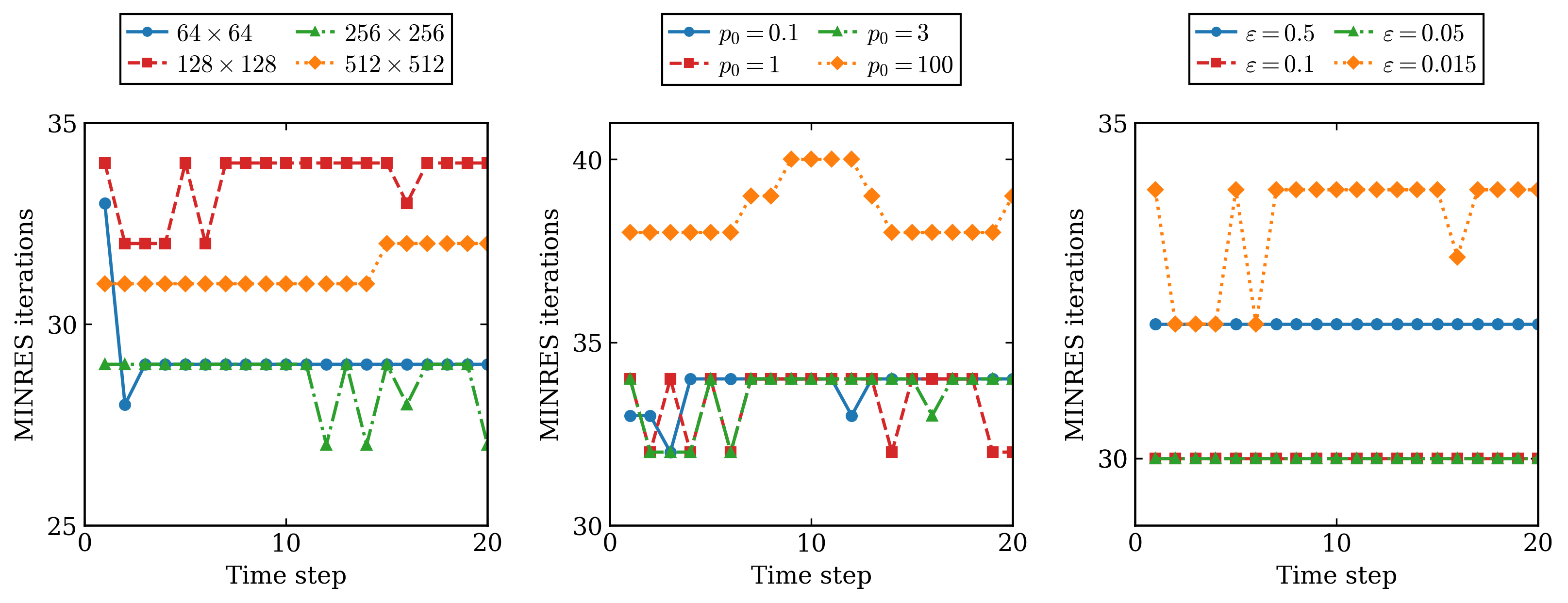}
    \caption{Ex. \ref{subsec:singel_tumor} Iteration counts over $20$ time steps for varying the grid points $n_{x_1},\, n_{x_2}$ (left), the proliferation growth parameter $p_0$ (middle), and the thickness parameter $\varepsilon$ (right). }
    \label{fig:HKNZ_iter_one}
\end{figure}

Figure~\ref{fig:HKNZ_iter_one} reports the \texttt{MINRES} iteration numbers over $20$ time steps for the HKNZ model with a single tumor growth simulations, varying the grid points $n_x,\, n_y$, the proliferation growth parameter $p_0$, and the thickness parameter $\varepsilon$. In all cases, the iteration counts remain bounded and as expected, they increase slightly at higher proliferation rates. Numerical experiments show a robust convergence behavior with respect to parameter changes.

\begin{figure}[htp!]
    \centering
    \includegraphics[width=0.88\linewidth]{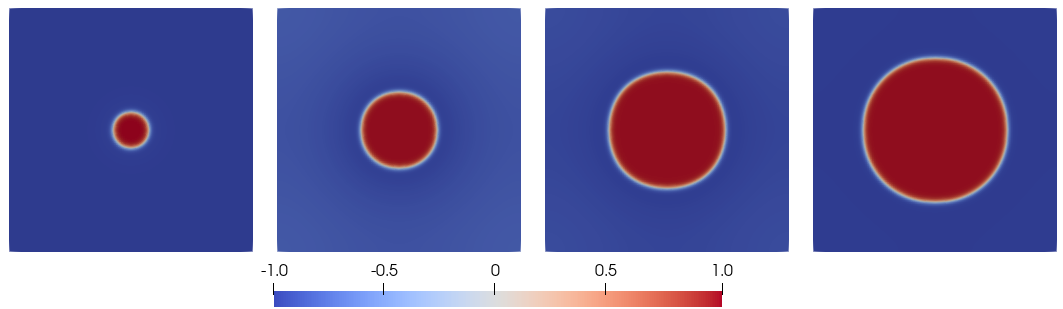}
     \includegraphics[width=0.88\linewidth]{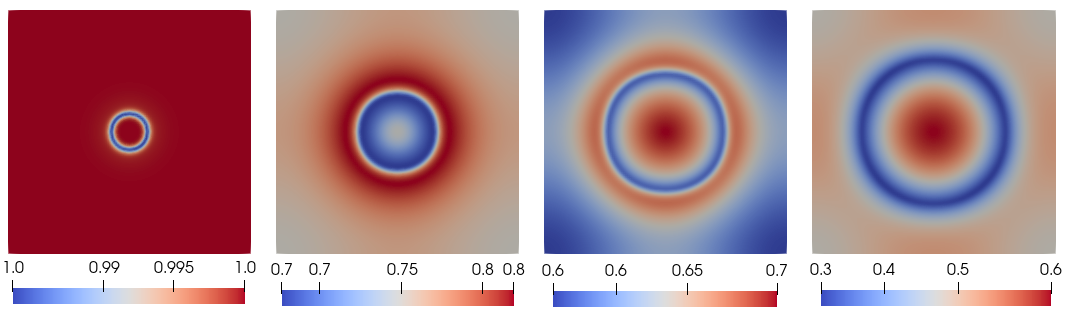}
    \caption{Ex. \ref{subsec:singel_tumor}: Tumor (top row) and nutrient (bottom row) concentration for time $t = 0.0001$, $t = 0.013$, $ t = 0.02$, $ t= 0.05$ from left to right. }
    \label{fig:one_tumor_u_sigma}
\end{figure}

Figure \ref{fig:one_tumor_u_sigma} shows how the nutrient is consumed by the tumor at times $t = 0.0001$, $0.013$, $0.02$, and $0.05$. In the nutrient plots, the colors represent different values for each plot.

\subsubsection{Three tumor growth simulation}\label{subsec:three_tumors}  
The initial condition consists of three tumors of distinct radii placed at different locations within the domain,
\begin{align}
  u_0(x_1, x_2)
  = & \tanh\!\left(\frac{0.15 - r_1}{\sqrt{2}\,\varepsilon}\right)+ \tanh\!\left(\frac{0.15 - r_2}{\sqrt{2}\,\varepsilon}\right) + \tanh\!\left(\frac{0.20 - r_3}{\sqrt{2}\,\varepsilon}\right)+ 2,
\end{align}
where 
\begin{align*}
  r_1 &= \sqrt{(x_1 - 0.45)^2 + (x_2 - 0.45)^2}, \;
  r_2 = \sqrt{(x_1 - 0.30)^2 + (x_2 + 0.50)^2}, \\
  r_3 &= \sqrt{(x_1 + 0.20)^2 + (x_2 + 0.20)^2}.
\end{align*}

\begin{figure}[htp!]
    \centering
    \includegraphics[width=0.88\linewidth]{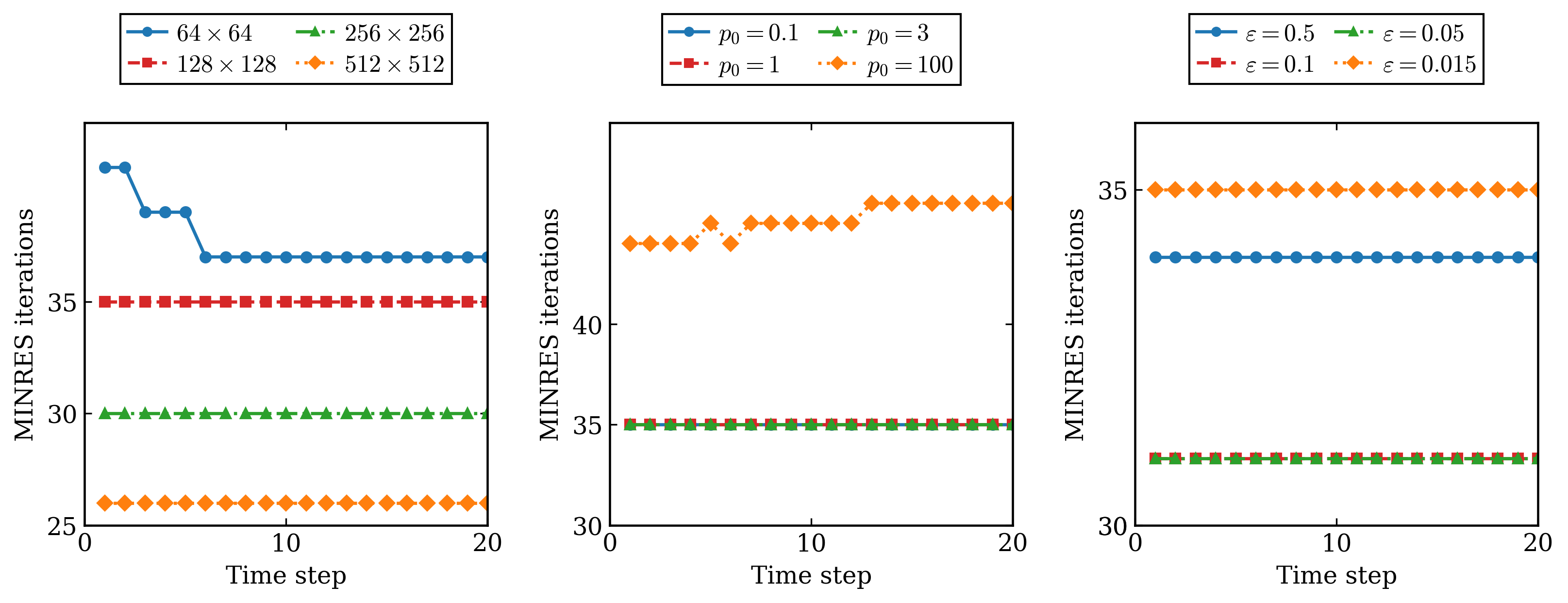}
    \caption{Ex. \ref{subsec:three_tumors} Iteration counts over $20$ time steps for varying the grid points $n_{x_1},\, n_{x_2}$ (left), the proliferation growth parameter $p_0$ (middle), and the thickness parameter $\varepsilon$ (right). }
    \label{fig:HKNZ_iter_three}
\end{figure}

In Figure~\ref{fig:HKNZ_iter_three}, we present the \texttt{MINRES} iteration counts over the first $20$ time steps for varying the grid resolution $(n_{x_1},n_{x_2})$, the proliferation growth parameter $p_0$, and the interface thickness parameter $\varepsilon$. As in Example~\ref{subsec:singel_tumor}, the iteration numbers remain in the bounded region, demonstrating the robustness of the proposed preconditioner. 

\begin{figure}[htp!]
    \centering
    \includegraphics[width=0.88\linewidth]{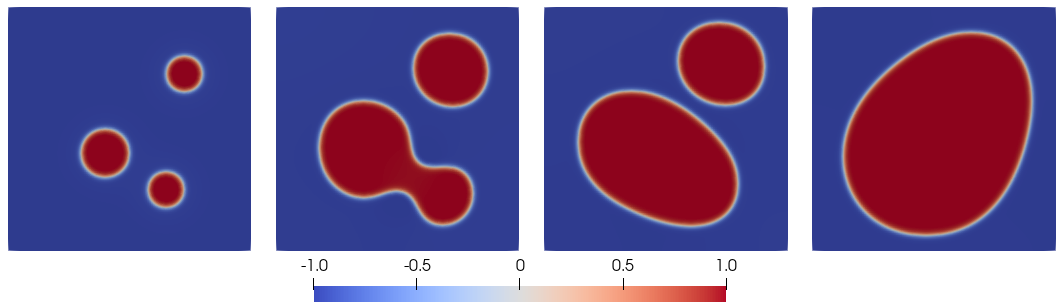}
     \includegraphics[width=0.88\linewidth]{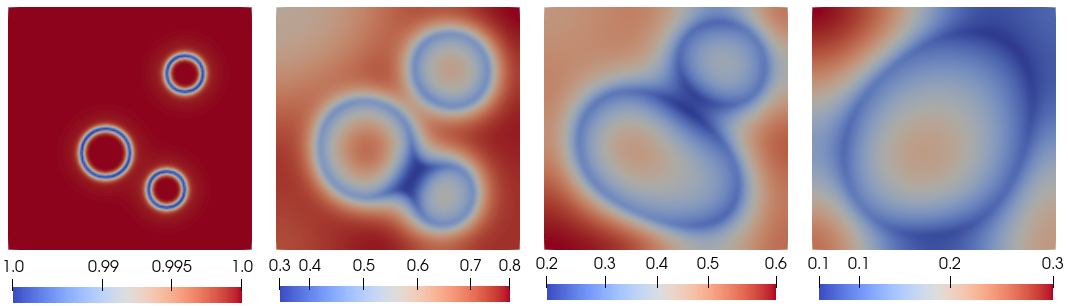}
    \caption{Ex. \ref{subsec:three_tumors} Three merging tumors (top row) and nutrient (bottom row) concentration for time $t = 0.0001$, $t = 0.03$, $ t = 0.07$, $ t= 0.2$ from left to right. }
    \label{fig:three_tumors}
\end{figure}

Figure~\ref{fig:three_tumors} illustrates the evolution of the three initially separated tumors. As time progresses, the tumors grow and eventually merge into a single larger tumor, while the nutrient concentration decreases inside the tumor region due to nutrient consumption and remains highest in the surrounding healthy tissue. These results are consistent with the expected qualitative behavior of the model.


\subsubsection{Two tumor growth simulation in 3D} \label{subsec:two_tumors}
\begin{figure}[htp!]
    \centering
    \includegraphics[width=0.88\linewidth]{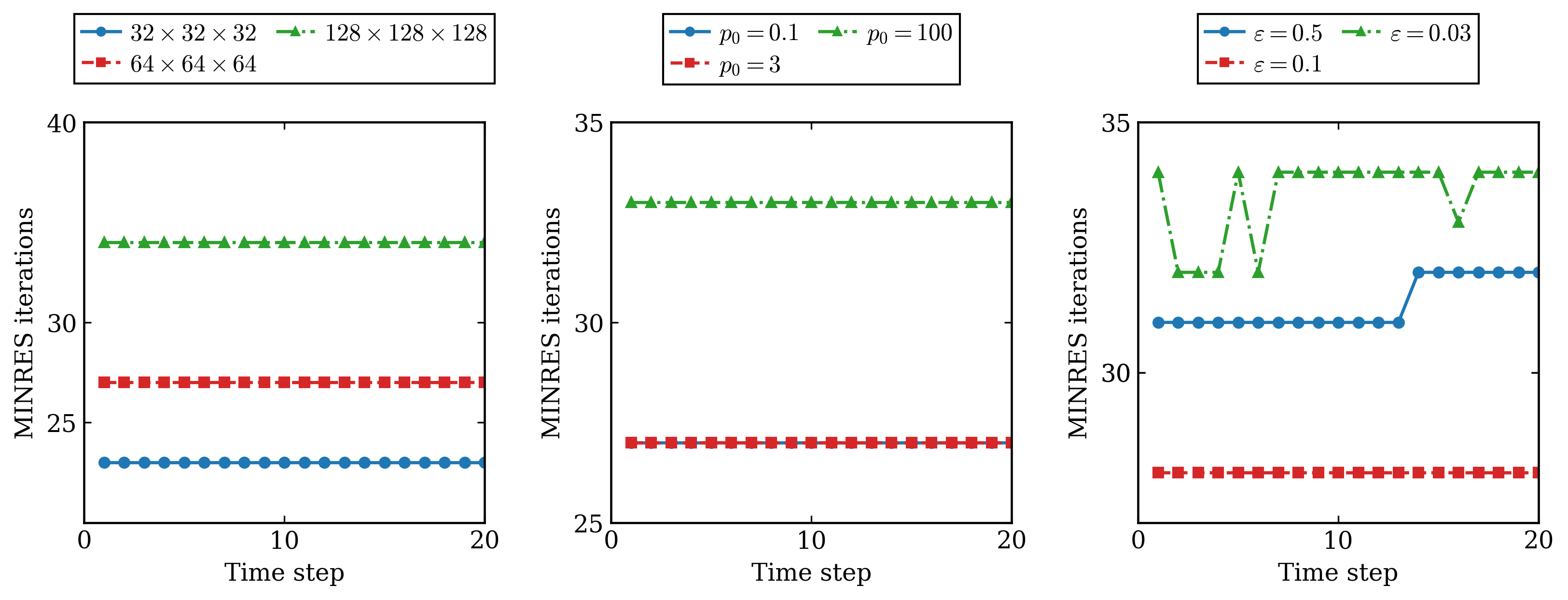}
    \caption{Ex. \ref{subsec:two_tumors} Iteration counts over $20$ time steps for varying the grid points $n_x,\, n_y\, n_z$ (left), the proliferation growth parameter $p_0$ (middle), and the thickness parameter $\varepsilon$ (right). }
    \label{fig:HKNZ_iter_two}
\end{figure}

The initial condition consists of two tumors of distinct radii, following the formulation of \cite{HikKNZ15}. In three dimensions, it reads
\begin{equation}\label{eq:two_tumor_ic}
  u_0(x_1, x_2, x_3)
  = 1
  - \tanh\!\left(\frac{r_1 - R_1}{\varepsilon\sqrt{2}}\right)
  - \tanh\!\left(\frac{r_2 - R_2}{\varepsilon\sqrt{2}}\right),
\end{equation}
where
\begin{align*}
  r_i &= \sqrt{(x_1 - c_1^{(i)})^2 + (x_2 - c_2^{(i)})^2 + x_3^2},
  \qquad i = 1, 2,
\end{align*}
and the tumor centers are given by
\begin{align*}
  c_1^{(1)} &= \phantom{-}\tfrac{1}{2}\sqrt{2}\!\left(R_1 + \tfrac{d}{2}\right), &
  c_2^{(1)} &= -\tfrac{1}{2}\sqrt{2}\!\left(R_1 + \tfrac{d}{2}\right), \\[4pt]
  c_1^{(2)} &= -\tfrac{1}{2}\sqrt{2}\!\left(R_2 + \tfrac{d}{2}\right), &
  c_2^{(2)} &= \phantom{-}\tfrac{1}{2}\sqrt{2}\!\left(R_2 + \tfrac{d}{2}\right).
\end{align*}
The parameters are set to $R_1 = 0.3$, $R_2 = 0.35$, and $d = 0.2$. Figure~\ref{fig:HKNZ_iter_two} shows the \texttt{MINRES} iteration counts over the first $20$ time steps for the three-dimensional two tumor example. As in the two-dimensional experiments, the proposed block preconditioner exhibits robust performance with respect to mesh refinement, the proliferation growth parameter $p_0$, and the interface thickness parameter $\varepsilon$. The iteration counts remain nearly constant throughout the simulation. These results demonstrate that the efficiency of the preconditioner is largely independent of both the discretization parameters and the physical model parameters.

\begin{figure}[htp!]
    \centering
     \includegraphics[width=0.19\linewidth]{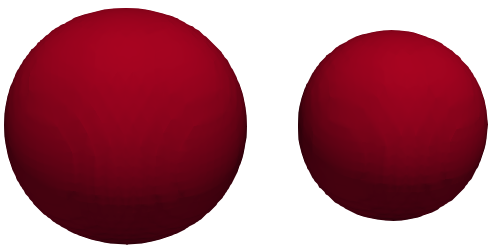}\hspace{0.4cm}
     \includegraphics[width=0.19\linewidth]{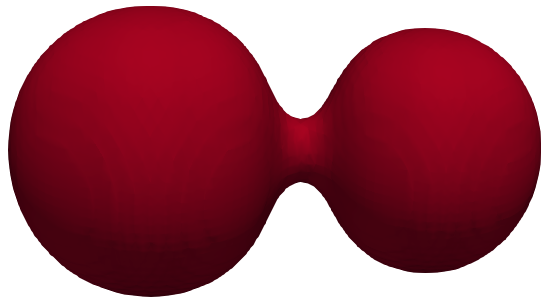}\hspace{0.4cm}
     \includegraphics[width=0.19\linewidth]{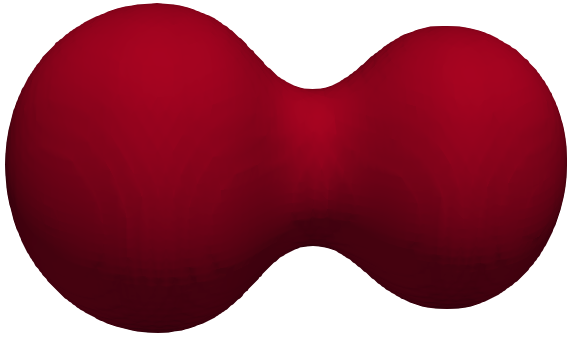}\hspace{0.4cm}
     \includegraphics[width=0.19\linewidth]{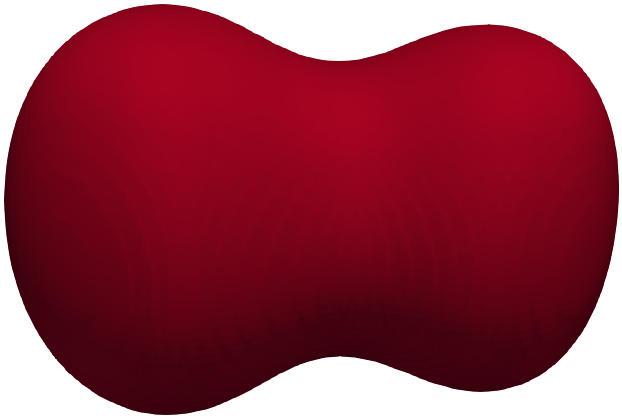}
    \vspace{0.25cm}
     
     \includegraphics[width=0.22\linewidth]{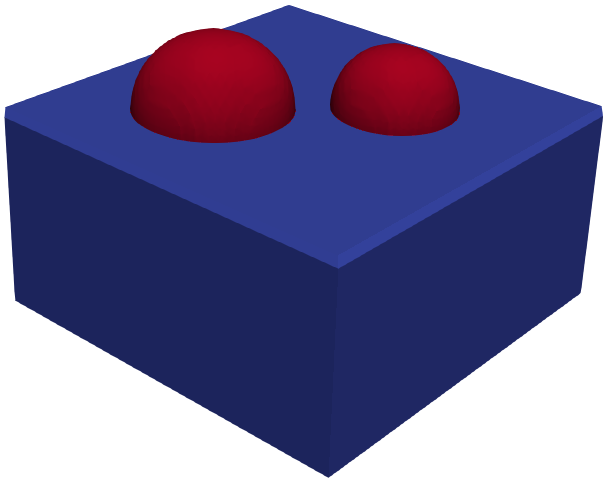}
     \includegraphics[width=0.22\linewidth]{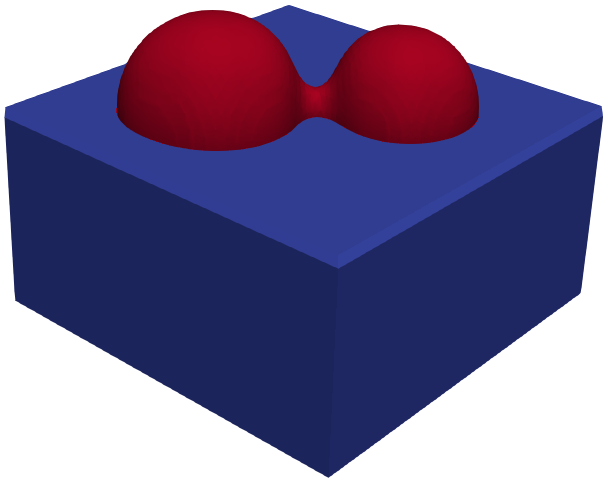}
     \includegraphics[width=0.22\linewidth]{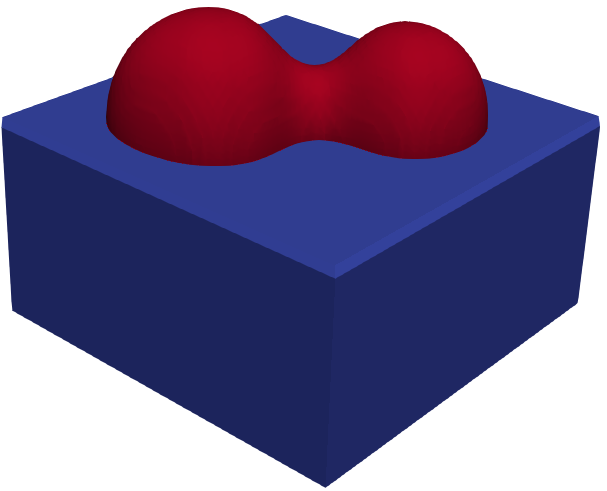}
     \includegraphics[width=0.22\linewidth]{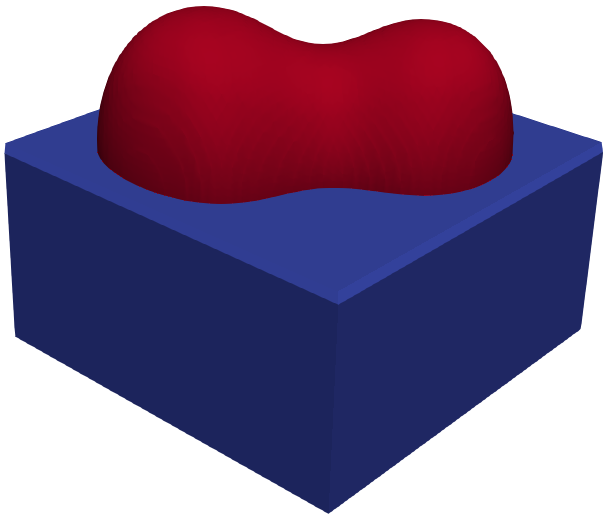}

     \includegraphics[width=0.22\linewidth]{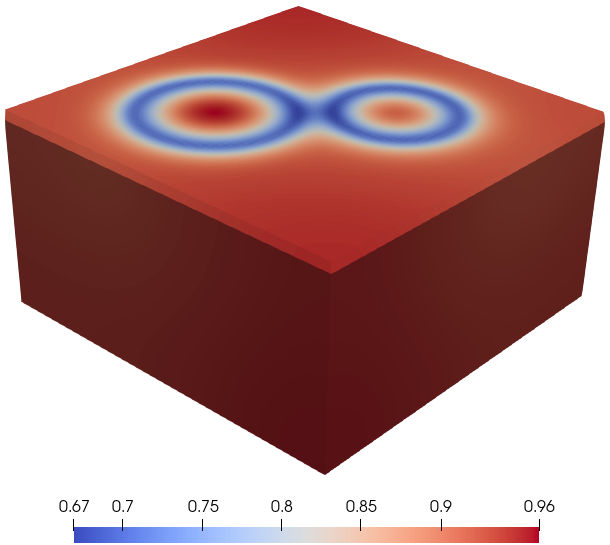}
     \includegraphics[width=0.22\linewidth]{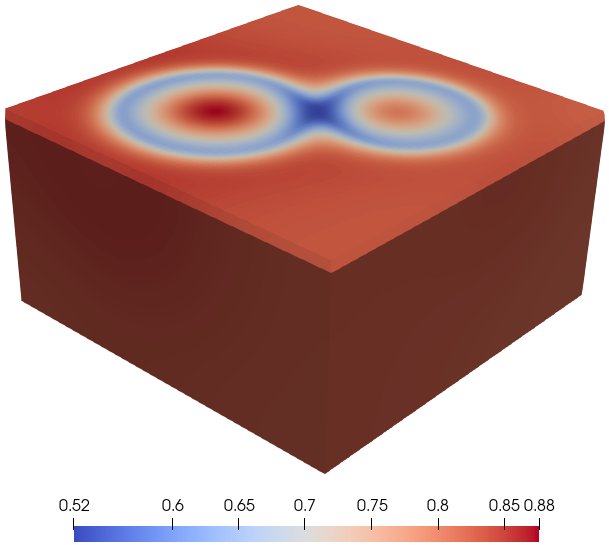}
     \includegraphics[width=0.22\linewidth]{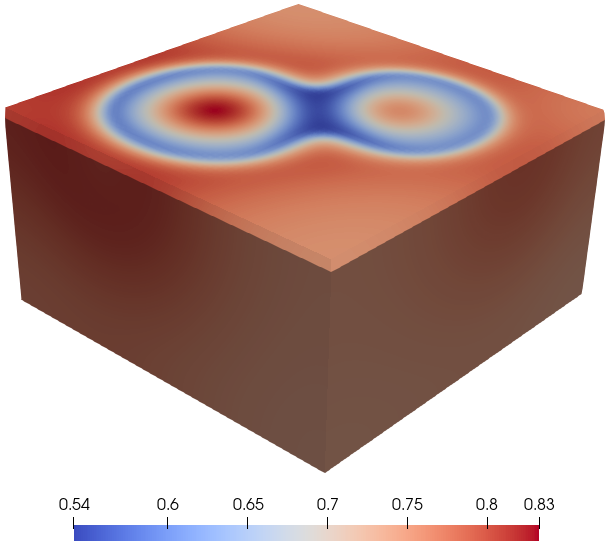}
     \includegraphics[width=0.22\linewidth]{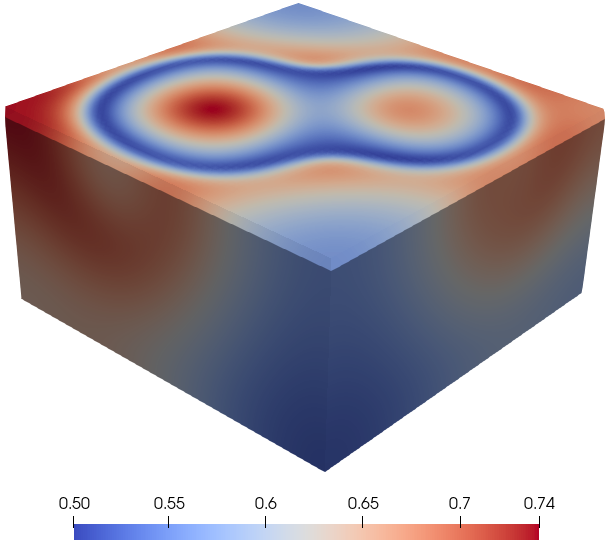}
    \caption{Ex. \ref{subsec:two_tumors} Two merging tumors (first and second rows) and nutrient (third row) concentration for time $t = 0.01$, $t = 0.0165$, $ t = 0.02$, $ t= 0.3$ from left to right. }
    \label{fig:two_tumors_3D}
\end{figure}

Figure~\ref{fig:two_tumors_3D} illustrates the evolution of the three-dimensional tumor configuration. Initially, two tumors of different sizes are separated by a small gap. As the simulation progresses, the tumors expand, interact, and eventually merge into a single connected tumor. As expected, the nutrient level is depleted inside the tumor due to consumption by proliferating cells and remains highest in the surrounding healthy tissue.

\begin{figure}[htp!]
    \centering
    \includegraphics[width=0.32\linewidth]{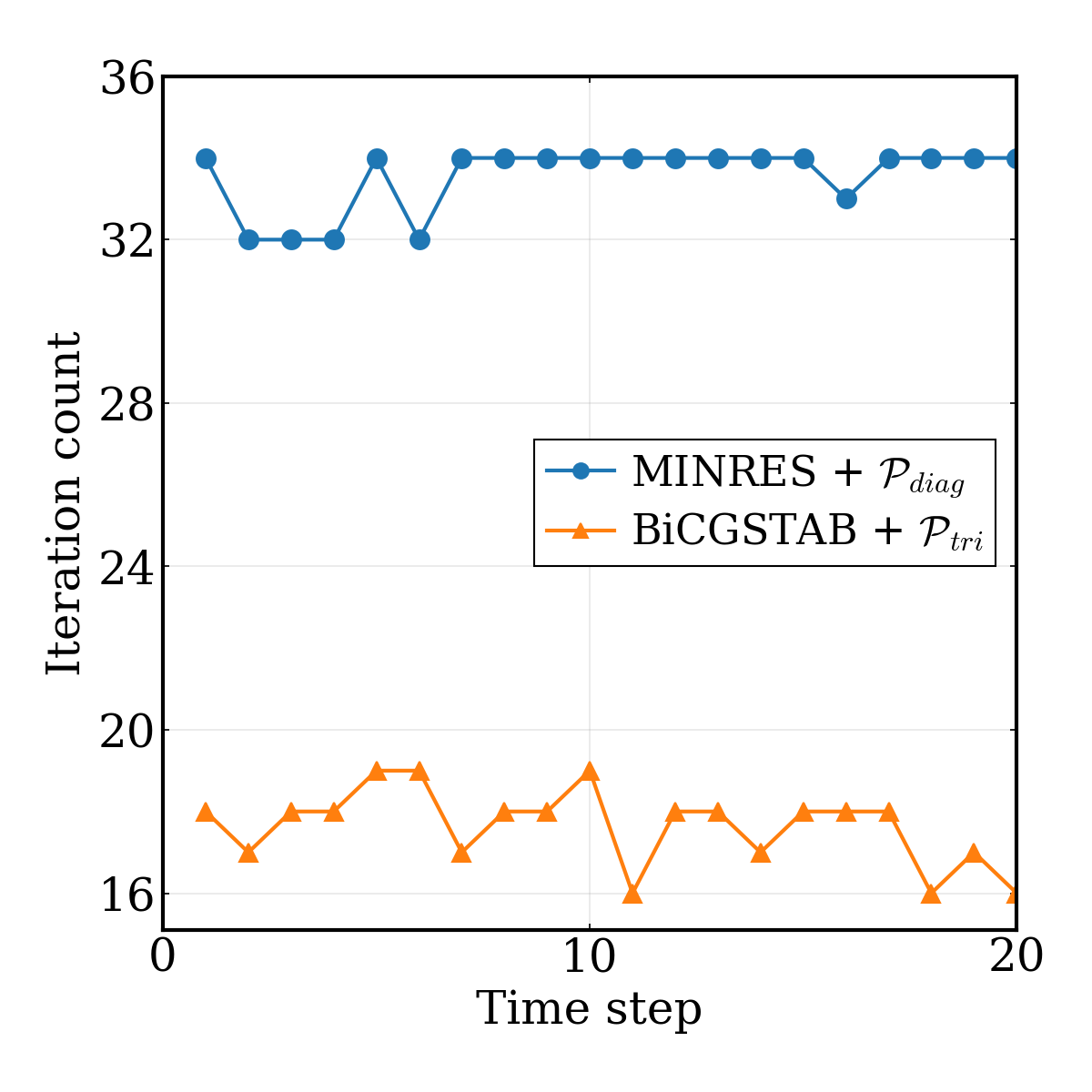}
    \includegraphics[width=0.32\linewidth]{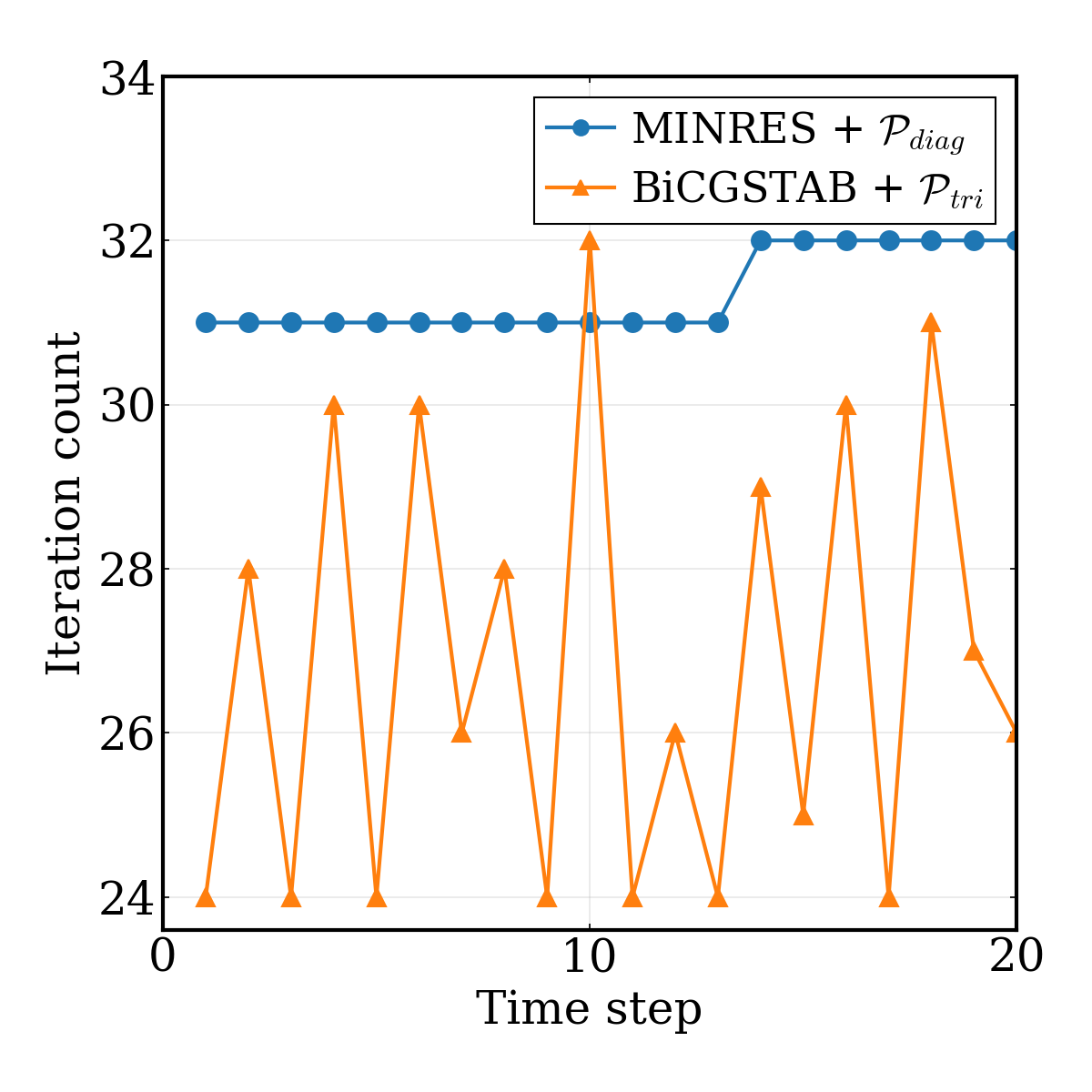}
    \caption{Number of Krylov iterations per time step over $20$ time steps for the smooth HKNZ model, comparing \texttt{MINRES} with the block-diagonal preconditioner $\bm{\mathcal{P}}_{diag}$ against \texttt{Bi-CGSTAB} with the block-triangular preconditioner $\bm{\mathcal{P}}_{tri}$: two-dimensional Ex. \ref{subsec:singel_tumor}  (left) and three-dimensional Ex. \ref{subsec:two_tumors}  (right).}
    \label{fig:HKNZ_iter_prec}
\end{figure}

Lastly, we also compare the performance of the block-diagonal preconditioner $\bm{\mathcal{P}}_{diag}$ against \texttt{Bi-CGSTAB} with the block-triangular preconditioner $\bm{\mathcal{P}}_{tri}$ in Figure~\ref{fig:HKNZ_iter_prec}. As expected, we see a modest improvements in iteration counts with the block-triangular preconditioner $\bm{\mathcal{P}}_{tri}$ in $2D$ and $3D$ cases.

\subsection{Model GLSS}
This section provides numerical results for the GLSS model and the smooth and nonsmooth variants differ the resulting linear systems that are solved using the block preconditioners described in Section~\ref{sec:prec}.

\subsubsection{Smooth GLSS}\label{subsec:garcke_smooth}
Here, we present numerical results for the smooth GLSS model on the domain $\Omega = (-12.5,12.5)^2$. To assess the robustness of the proposed method with respect to different interface geometries, we consider three representative initial conditions: a four-finger shape, a superellipse, and a Cassini oval. The corresponding phase-field functions are defined as follows:

\begin{itemize}
\item[(a)] Superellipse:
\begin{equation*}
\begin{split}
    u_0(x_1,x_2) &= \tanh\!\left( \frac{R_0-r(x_1,x_2)}
    {\sqrt{2}\,\varepsilon} \right), \;
    r(x_1,x_2)  = \left(|x_1|^{5/2}+|x_2|^{5/2}\right)^{2/5}, 
\end{split}
\end{equation*}
\item[(b)] Four-finger shape:
\begin{equation*}
\begin{split}
    u_0(x_1,x_2) &= \tanh\!\left( \frac{R(\theta)-r}{\sqrt{2}\,\varepsilon} \right), \;
    R(\theta)=R_0\left(1+0.4\cos(4\theta)+0.1\cos(8\theta)\right),  \\
    r(x_1,x_2) & =\sqrt{x_1^2+x_2^2}, \quad \theta=\arctan(x_2,x_1)-\pi/4,
\end{split}
\end{equation*}
\item[(c)] Cassini oval:
\begin{equation*}
\begin{split}
    u_0(x_1,x_2) &= \tanh\!\left( -\frac{C(x_1,x_2)} {b^2\sqrt{2}\,\varepsilon} \right), \; 
    C(x_1,x_2) = \left(x_1^2+x_2^2+a^2\right)^2 - 4a^2x_1^2-b^4, 
\end{split}
\end{equation*}
\end{itemize}
where $R_0 = 2.0$, $a=2.0$, and $b=2.05$. Unless otherwise stated, the numerical experiments presented in this section are performed using the superellipse initial condition. The initial nutrient concentration is taken as
$\sigma_0(x_1,x_2)=1$.

\begin{figure}[htp!]
    \centering
    \includegraphics[width=0.88\linewidth]{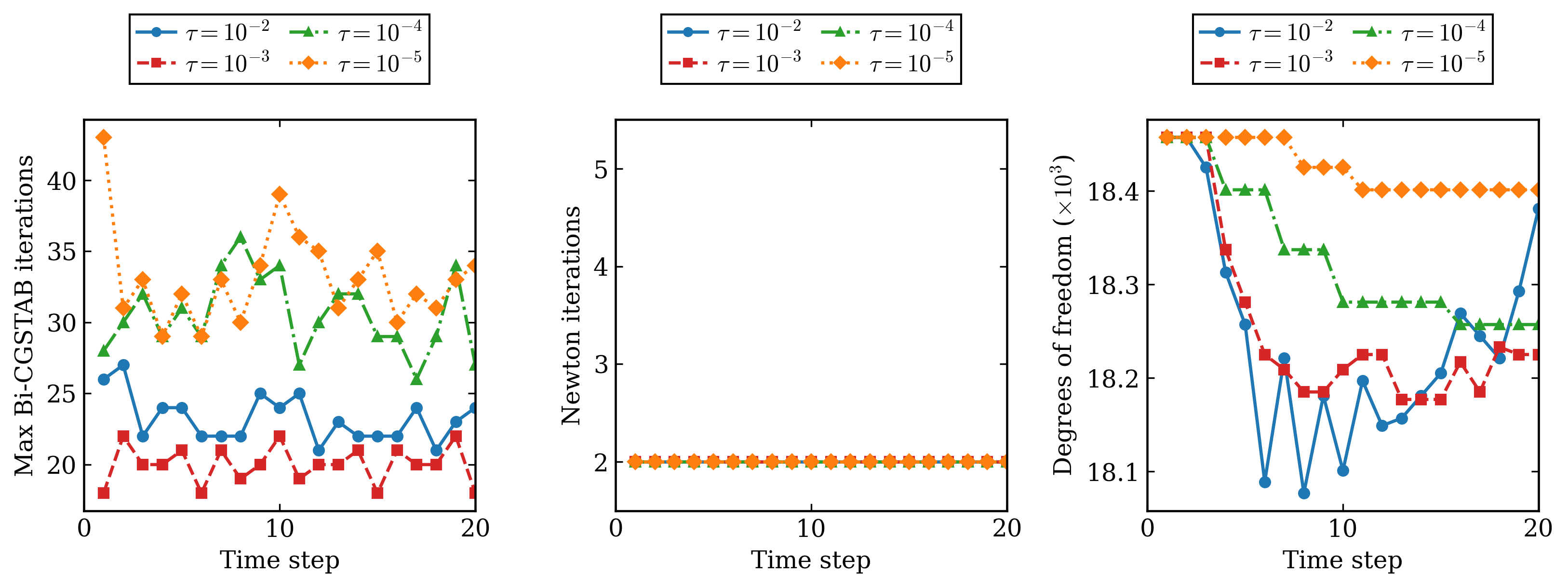}
     \includegraphics[width=0.88\linewidth]{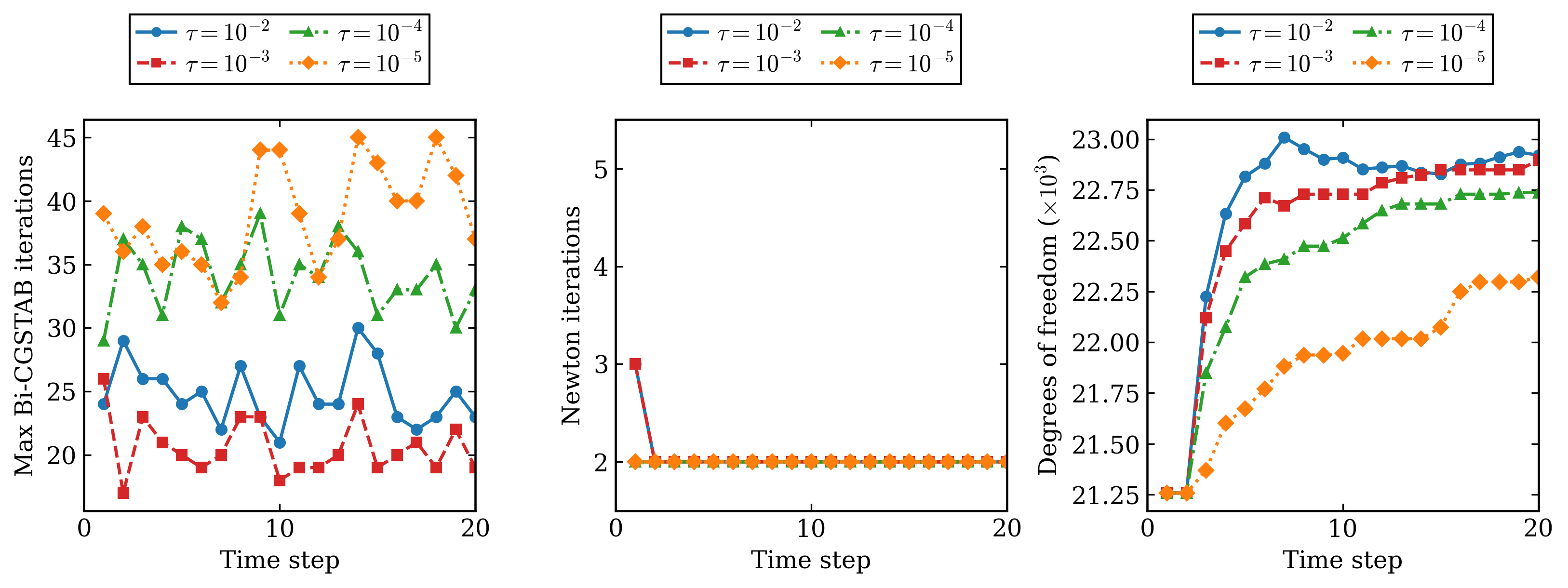}
      \includegraphics[width=0.88\linewidth]{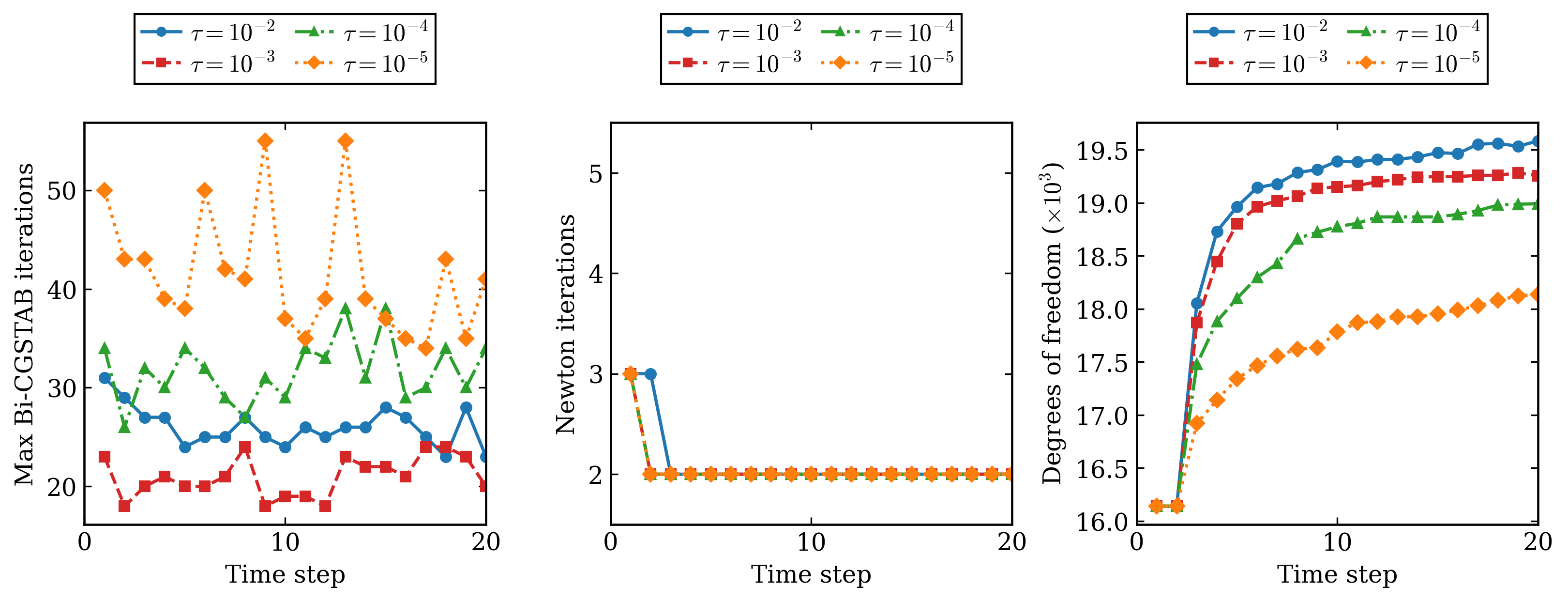}
    \caption{Ex.~\ref{subsec:garcke_smooth} Iteration counts over $20$ time steps for varying time-step size $\tau$: maximum \texttt{Bi-CGSTAB} iterations (left), Newton iterations (center), and degrees of freedom (right) with the initial condition (a), (b), (c) from top to bottom. The linear and nonlinear solver tolerance are set to  $\texttt{rtol} = 10^{-12}$ for \texttt{Bi-CGSTAB} and $\texttt{rtol} = 10^{-10}$ and $\texttt{rtol} = 10^{-5}$ for Newton.}
    \label{fig:GLSS_iter_a_new}
\end{figure}

We fix the coarse and fine mesh sizes $h_c = 25\cdot2^{-6}$ and $h_f = 25\cdot2^{-10}$, the refinement threshold $\xi = 0.075$. All parameters used in the simulations are as given in Table~\ref{tab:garcke_params}. We compare four time-step sizes $\tau\in\{10^{-2},10^{-3},10^{-4},10^{-5}\}$ and report, over the first 20 time steps, the maximum \texttt{Bi-CGSTAB} iteration count per time step, the Newton iteration count, and the number of degrees of freedom obtained by adaptive mesh refinement in Figure~\ref{fig:GLSS_iter_a_new} with the initial condition (a), (b), and (c). 
\begin{table}[htp!]
\centering
\caption{Parameter values used for the smooth GLSS model experiments.}
\label{tab:garcke_params}
\begin{tabular}{ccccc}
\hline
$\tau = 10^{-3} $ & $\beta = 0.1 $ & $\varepsilon= 0.02$ & $\lambda= 0.03$ & $\mathcal{C} = 2.0$ \\
$D=1.0$ & $\chi = 5.0$ & $\mathcal{P}=0.5$ & $\mathcal{A} = 0 $ &  \\
\hline
\end{tabular}
\end{table}

\begin{figure}[htp!]
    \centering
    \includegraphics[width=0.88\linewidth]{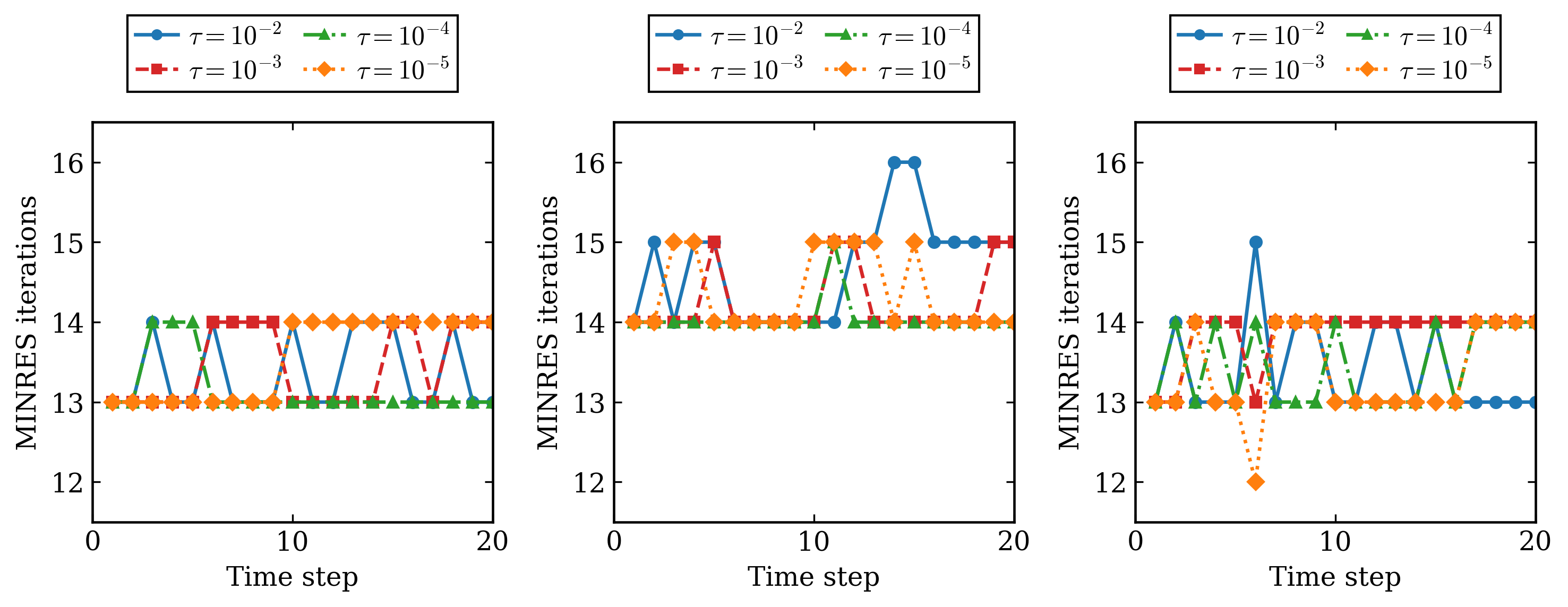}
    \caption{Ex.~\ref{subsec:garcke_smooth} Iteration counts over $20$ time steps for varying time-step size $\tau$: \texttt{MINRES} iterations with the initial condition (a), (b), (c) from left to right for the solution of the system \eqref{GLSS_nonsmooth_matrix_2}.}
    \label{fig:GLSS_iter_sigma}
\end{figure}

Figure~\ref{fig:GLSS_iter_a_new} reports the maximum \texttt{Bi-CGSTAB} iterations, Newton iterations, and degrees of freedom over the first $20$ time steps for different time-step sizes $\tau$ and the three initial conditions. 
The tolerances are specified in the caption. We have numerically observed that the iteration counts are reduced compared to tighter tolerances, without affecting the qualitative behavior of the problem. The Newton method converges in only two iterations throughout almost the entire simulation, demonstrating the robustness of the nonlinear solver. Although the number of \texttt{Bi-CGSTAB} iterations increases as the time-step size decreases, it remains within a relatively small range. The adaptive mesh evolves differently for each value of $\tau$, resulting in different numbers of degrees of freedom.  However,  this has only a limited impact on the performance of the solver. Figure~\ref{fig:GLSS_iter_sigma} demonstrates that the number of \texttt{MINRES} iterations is nearly independent of the time-step size for the solution of the system \eqref{GLSS_nonsmooth_matrix_2}. Overall, these results indicate that the proposed preconditioners remain effective over a wide range of time-step sizes and for all three initial conditions.

\begin{figure}[htp!]
    \centering
    \includegraphics[width=0.88\linewidth]{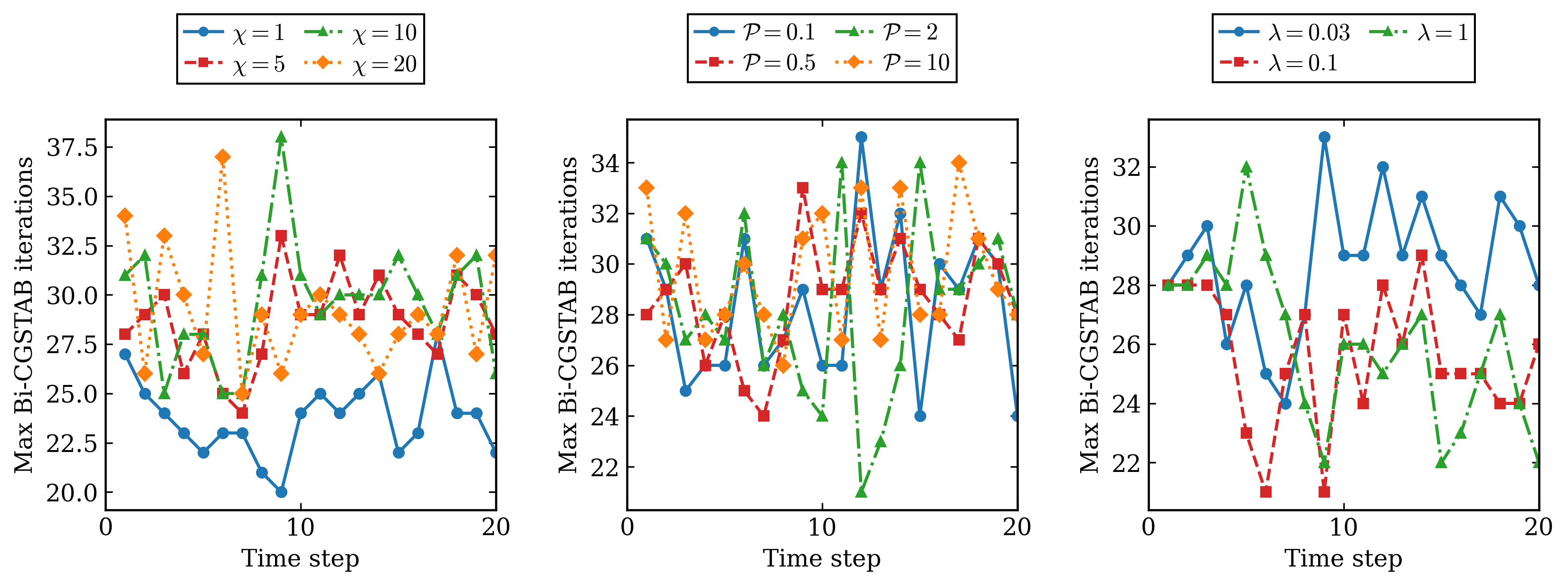}
    \caption{Ex.~\ref{subsec:garcke_smooth} Maximum \texttt{Bi-CGSTAB} iterations $20$ time steps for varying $\chi$ (left), $\mathcal{P}$ (middle), and $\lambda$ (right) with the initial condition (b).}
    \label{fig:GLSS_iter_bicgstab_parameter_b}
\end{figure}

Figure~\ref{fig:GLSS_iter_bicgstab_parameter_b} shows the maximum \texttt{Bi-CGSTAB} iteration counts for different values of the chemotaxis parameter $\chi$, the proliferation rate $\mathcal{P}$, and the nutrient consumption parameter $\lambda$ using the initial condition~(b). The iteration counts remain within a relatively small range over all time steps, indicating the robustness of the proposed preconditioner with respect to variations in the physical model parameters.

\begin{figure}[htp!]
    \centering
    \includegraphics[width=0.88\linewidth]{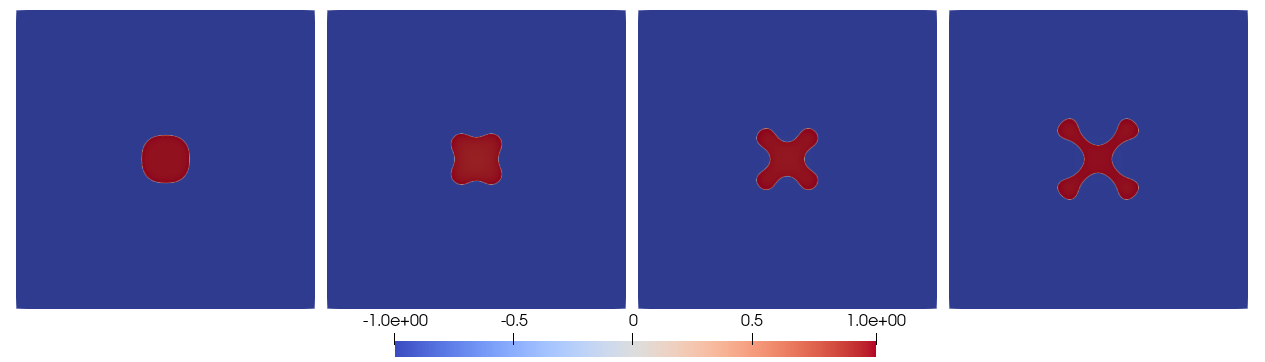}
    \includegraphics[width=0.88\linewidth]{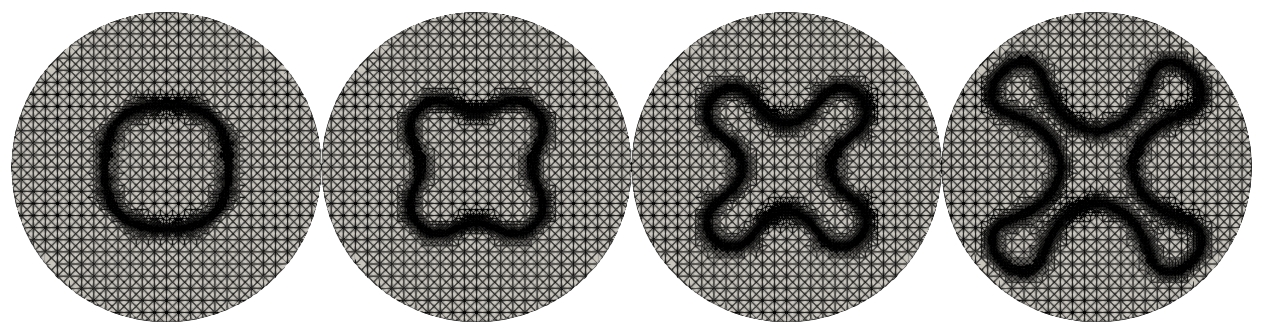}
    \includegraphics[width=0.88\linewidth]{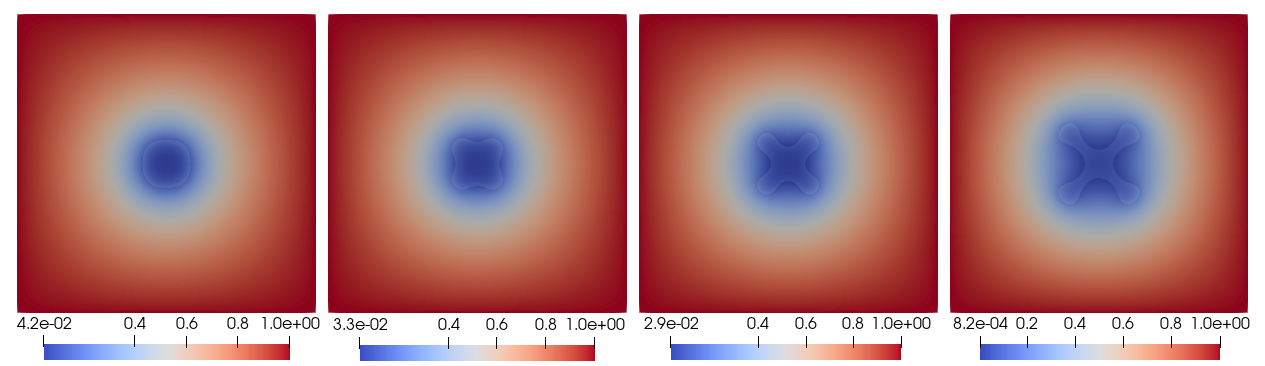}
    \caption{Ex.~\ref{subsec:garcke_smooth} Tumor (first), mesh (second rows) and nutrient (third row) concentration for time $t = 0.01$, $t = 3$, $ t = 5$, $ t= 8$ from left to right with the initial condition (a). }
    \label{fig:GLSS_u_sigma_a}
\end{figure}

\begin{figure}[htp!]
   \centering
   \includegraphics[width=0.88\linewidth]{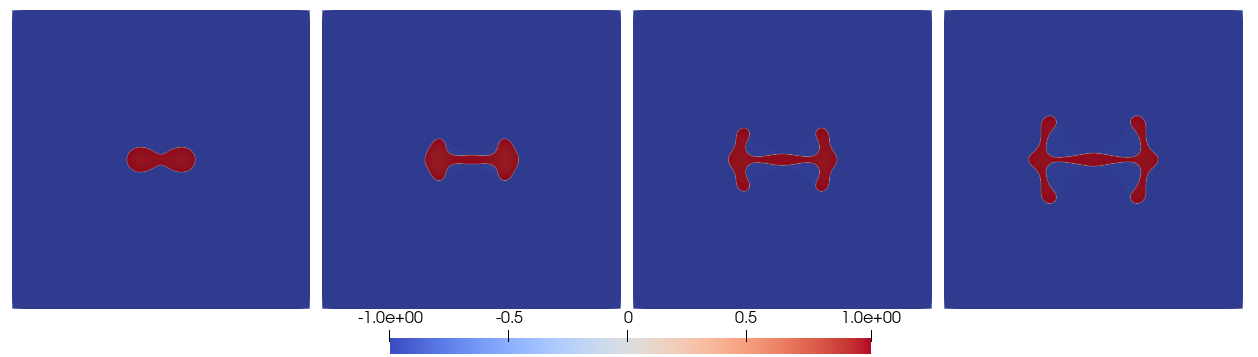}
   \includegraphics[width=0.88\linewidth]{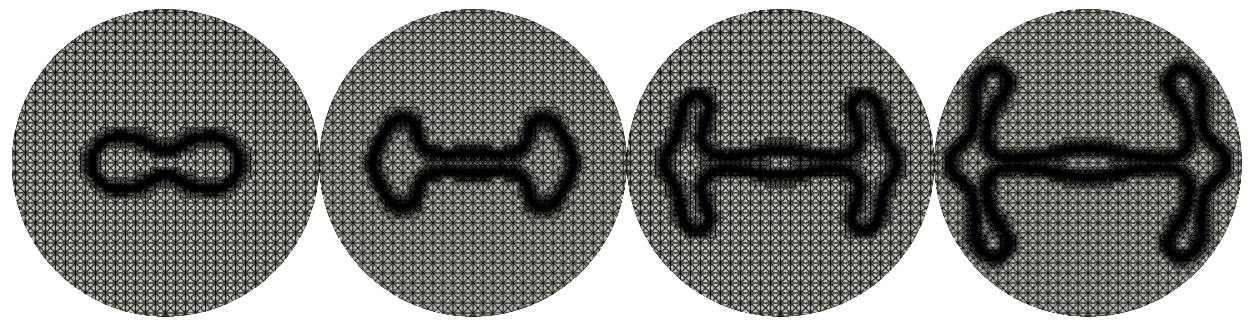}
   \includegraphics[width=0.88\linewidth]{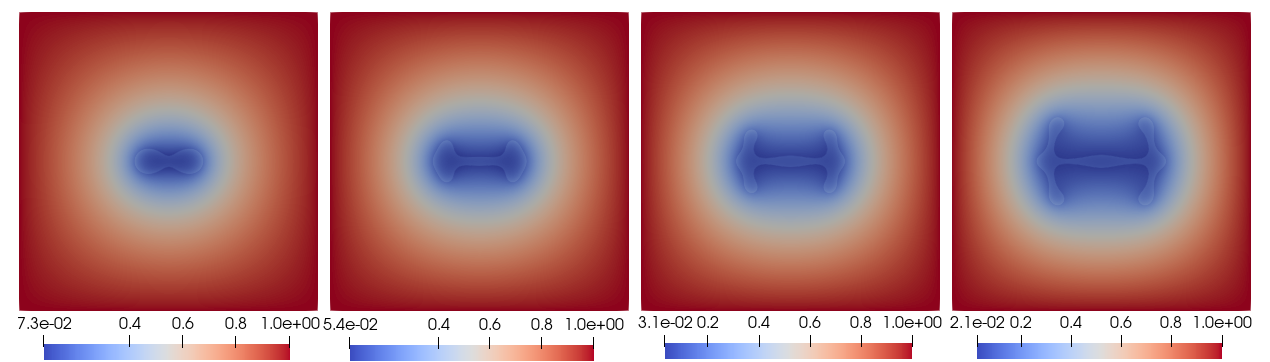}
   \caption{Ex.~\ref{subsec:garcke_smooth} Tumor (first), mesh (second rows) and nutrient (third row) concentration for time $t = 0.01$, $t = 5$, $ t = 8$, $ t= 11$ from left to right with the initial condition (c). }
   \label{fig:GLSS_u_sigma_c}
\end{figure}

Figures~\ref{fig:GLSS_u_sigma_a} and~\ref{fig:GLSS_u_sigma_c}  shows the evolution of the tumor, the adaptive mesh, and the nutrient concentration for  initial conditions (a) and (c). As the tumor grows, it develops fingers towards regions with higher nutrient concentration, allowing it to access additional nutrients while maintaining its overall symmetry. Throughout the simulation, the adaptive mesh remains concentrated around the diffuse interface, providing higher resolution only where it is needed. We report that qualitatively similar results are obtained for the other initial condition (b).

We also extend the numerical experiments of the smooth GLSS model to three spatial dimensions. The computational domain is taken as $\Omega = (-3, 3)^3$, and  all parameters are as given in Table~\ref{tab:garcke_params}, except the chemotaxis parameter $\chi=30$ which causes quicker tumor evolutions. The coarse mesh size is $h_c = 6 \cdot (1/14) $ and the fine mesh size is $h_f = 6 \cdot (1/14) \cdot 2^{-4}$. 

As initial condition, we use a three-dimensional Cassini oval obtained by
revolving the two-dimensional Cassini curve around the $x_1$-axis,
\begin{equation*}
  u_0(x_1, x_2, x_3)
  = \tanh\!\left( - \frac{C(x_1, r)}{\,b^2\sqrt{2}\,\varepsilon} \right), \qquad
  C(x_1, r) = \bigl(x_1^2 + r^2 + a^2\bigr)^2 - 4a^2 x_1^2 - b^4,
\end{equation*}
where $r = \sqrt{x_2^2 + x_3^2}$ denotes the distance from the $x_1$-axis, and the parameters $a = 2.0$ and $b = 2.05$.

\begin{figure}[htp!]
    \centering
    \includegraphics[width=0.9\linewidth]{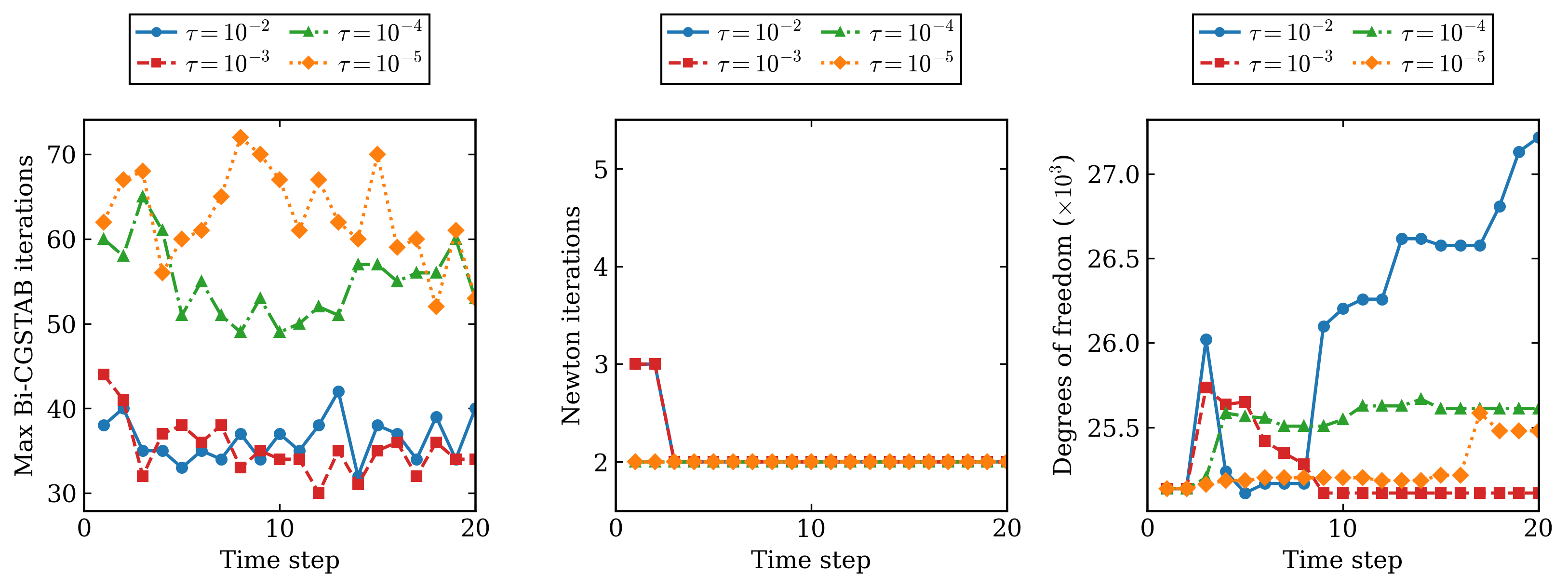}
    \caption{Ex.~\ref{subsec:garcke_smooth} Iteration counts over $20$ time steps for varying time-step size $\tau$: maximum \texttt{Bi-CGSTAB} iterations (left), Newton iterations (center), and degrees of freedom (right) with the initial condition (c) in 3D with \texttt{AMG} preconditioning.}
    \label{fig:GLSS_iter_c_3D_AMG}
\end{figure}

\begin{figure}[htp!]
    \centering
    \includegraphics[width=0.9\linewidth]{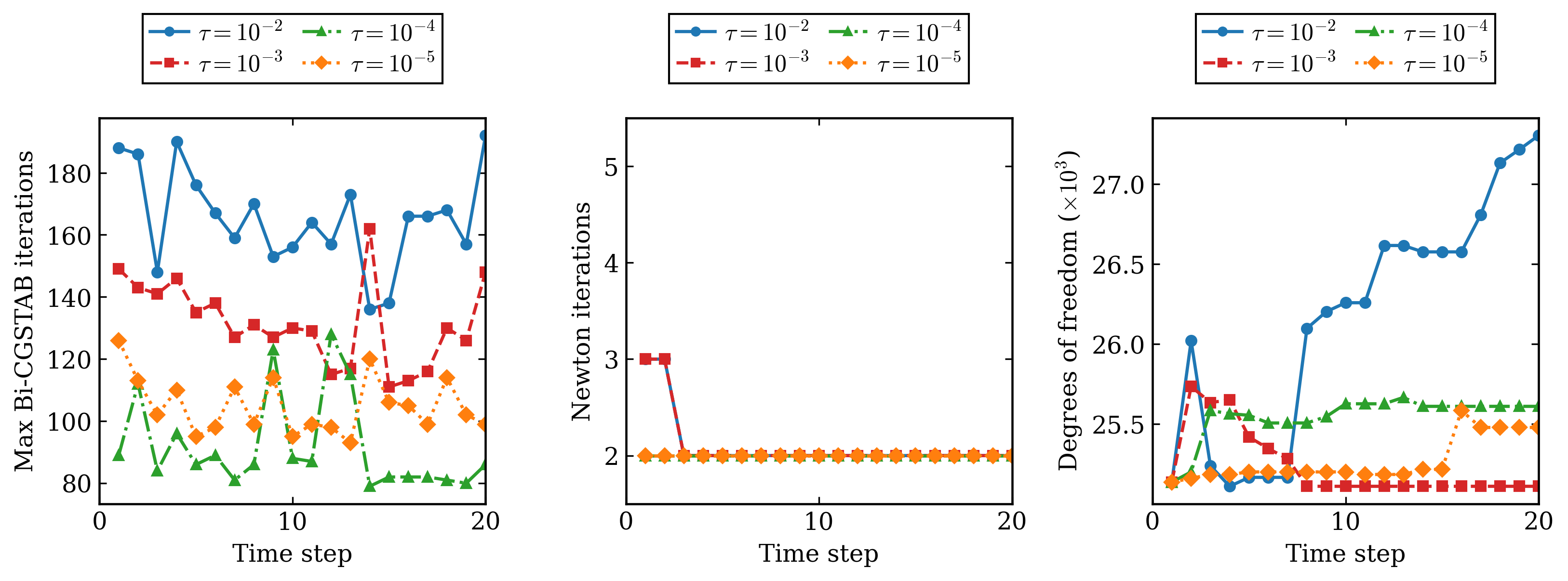}
    \caption{Ex.~\ref{subsec:garcke_smooth} Iteration counts over $20$ time steps for varying time-step size $\tau$: maximum \texttt{Bi-CGSTAB} iterations (left), Newton iterations (center), and degrees of freedom (right) with the initial condition (c) in 3D with \texttt{ILU} preconditioning.}
    \label{fig:GLSS_iter_c_3D_ILU}
\end{figure}

\begin{figure}[htp!]
    \centering
     \includegraphics[width=0.19\linewidth]{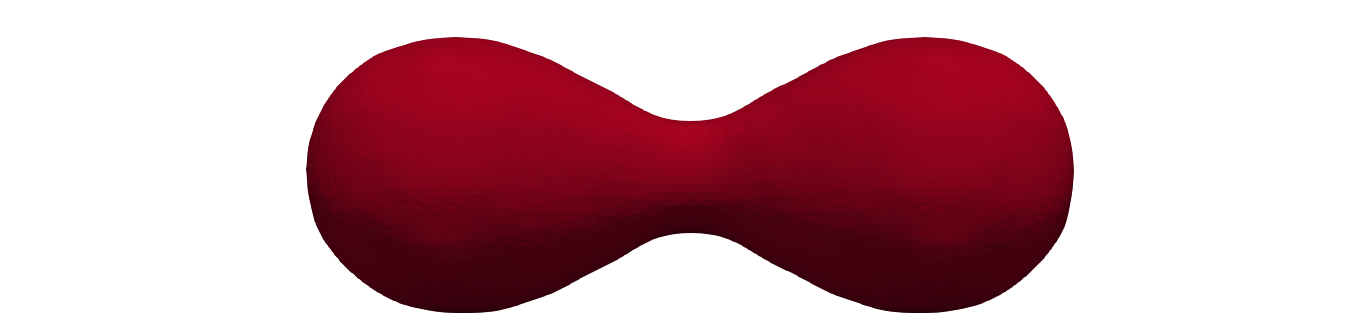}\hspace{0.4cm}
     \includegraphics[width=0.19\linewidth]{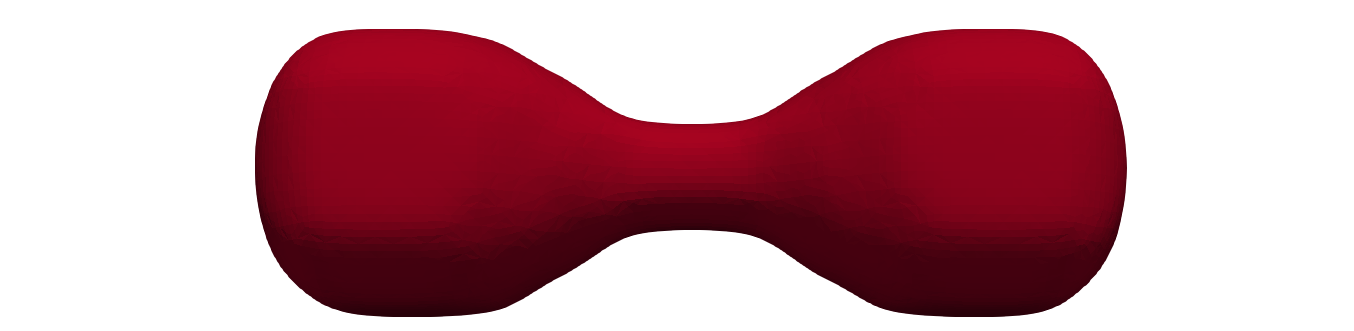}\hspace{0.4cm}
     \includegraphics[width=0.19\linewidth]{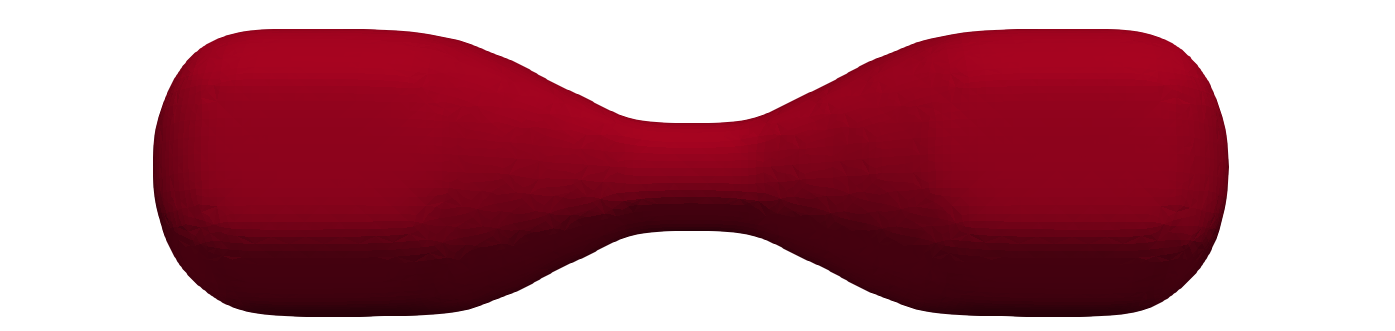}\hspace{0.4cm}
     \includegraphics[width=0.19\linewidth]{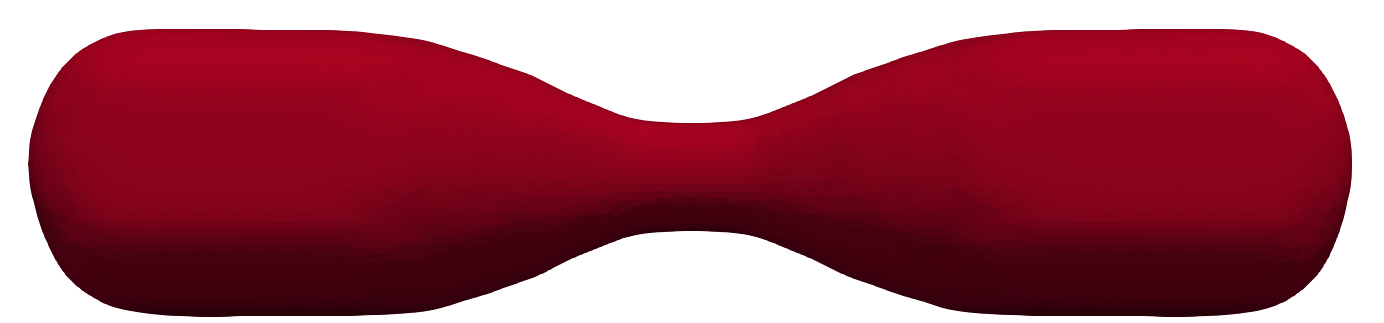}
    \vspace{0.25cm}
     
     \includegraphics[width=0.22\linewidth]{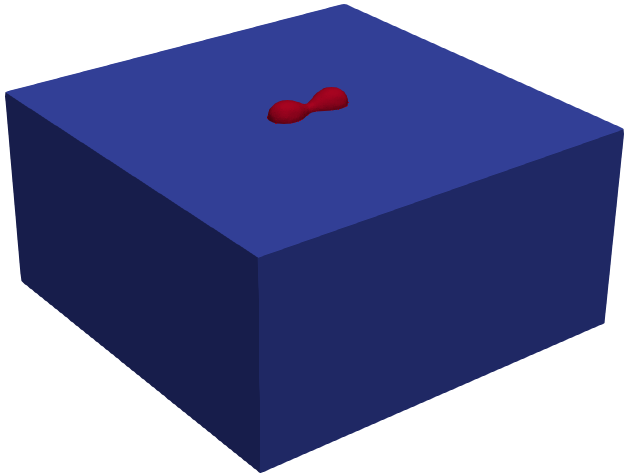}
     \includegraphics[width=0.22\linewidth]{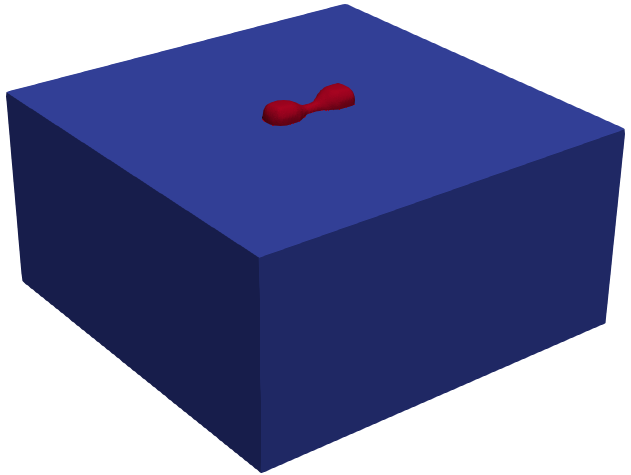}
     \includegraphics[width=0.22\linewidth]{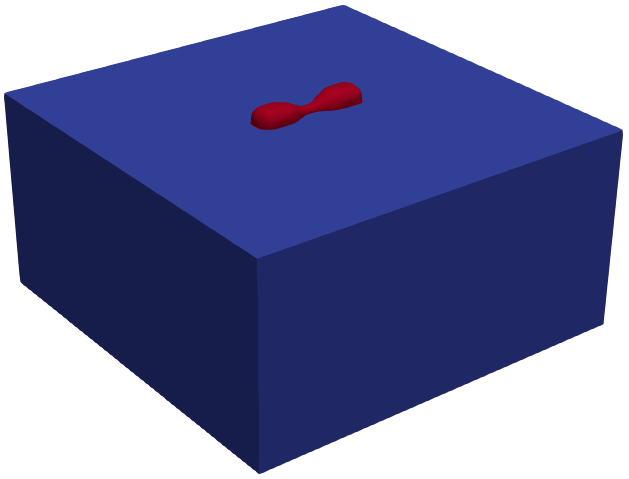}
     \includegraphics[width=0.22\linewidth]{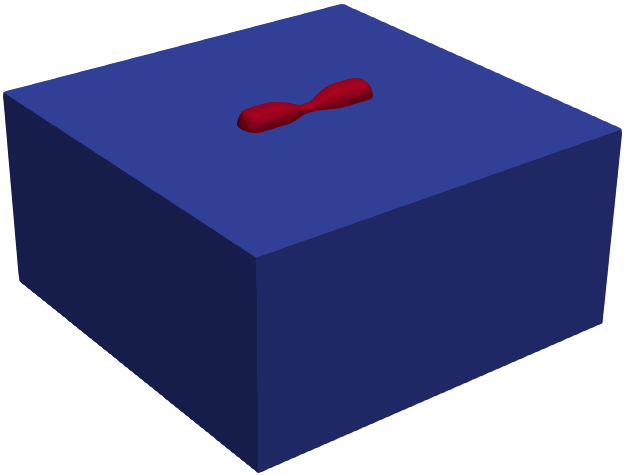}

     \includegraphics[width=0.22\linewidth]{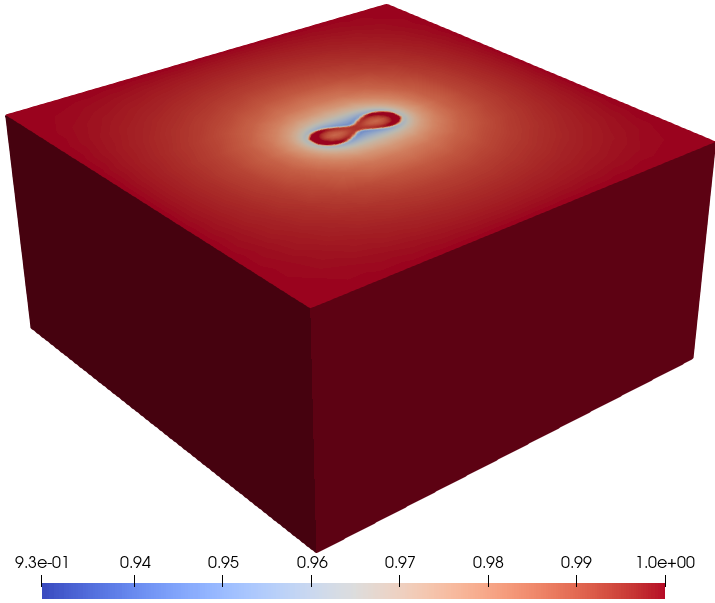}
     \includegraphics[width=0.22\linewidth]{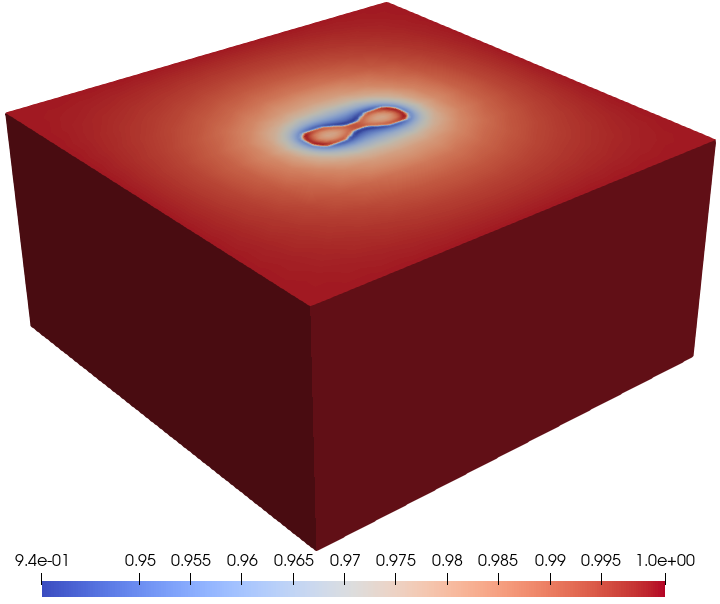}
     \includegraphics[width=0.22\linewidth]{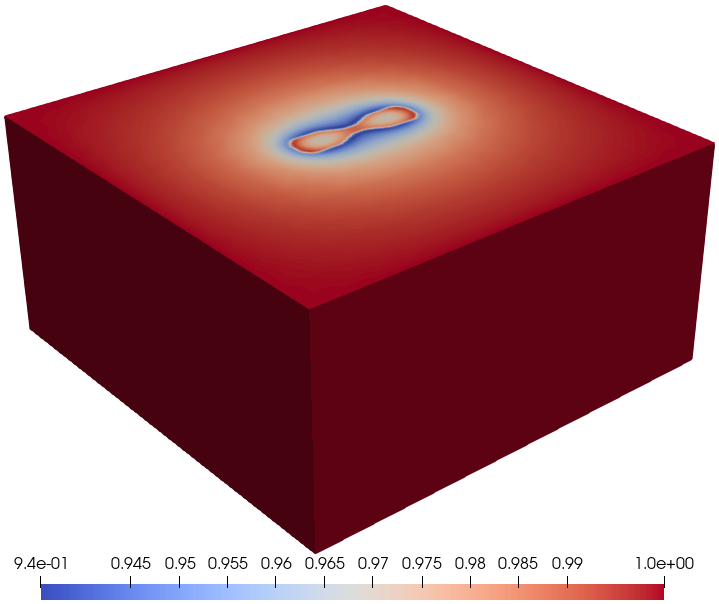}
     \includegraphics[width=0.22\linewidth]{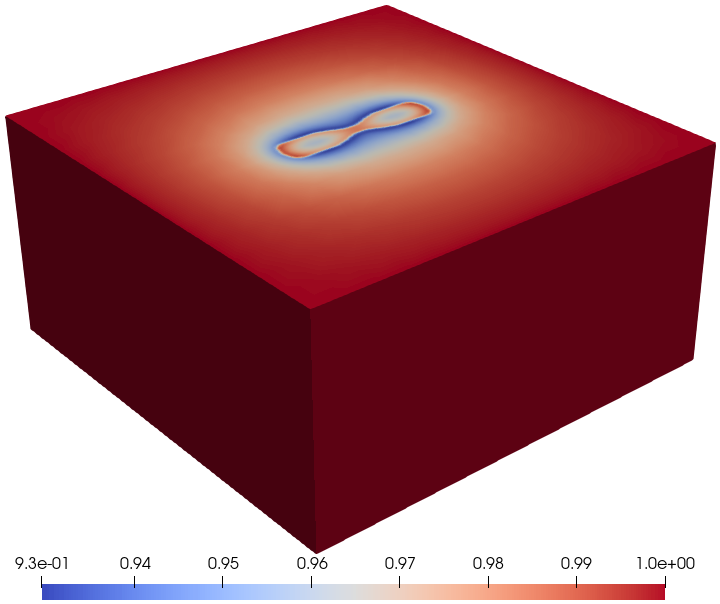}
    \caption{Ex. \ref{subsec:garcke_smooth} Tumor (first and second rows) and nutrient (third row) concentration for time $t = 0.01$, $t = 0.5$, $ t = 1$, $ t= 1.5$ from left to right. }
    \label{fig:Garkce_tumors_3D}
\end{figure}

In three dimensions, we compare \texttt{ILU} and \texttt{AMG} as approximate solvers for the diagonal blocks of the preconditioner. We therefore replace \texttt{ILU} by a single V-cycle of smoothed aggregation \texttt{AMG}, as implemented in \texttt{PyAMG}~\cite{pyamg2023}, using classical strength-of-connection with threshold $\theta = 0.25$ and symmetric Gauss--Seidel pre- and post-smoothing. The \texttt{AMG} hierarchy is constructed with a maximum of $20$ levels and a coarsest-level size of at most $50$ unknowns.
Figures~\ref{fig:GLSS_iter_c_3D_AMG} and~\ref{fig:GLSS_iter_c_3D_ILU}  compare the solver performance of \texttt{AMG} and \texttt{ILU} preconditioning for the three-dimensional smooth GLSS model. The Newton iteration counts and degree-of-freedom histories are identical in both cases, confirming that the choice of preconditioner does not affect the nonlinear convergence or the adaptive mesh refinement. The \texttt{Bi-CGSTAB} iteration counts differ substantially between the two  preconditioners: \texttt{AMG} requires roughly $30$--$85$ iterations, while \texttt{ILU} requires $80$--$220$, with considerably more variation across time steps. This confirms that \texttt{AMG} is the more effective choice for three-dimensional problems.  Figure~\ref{fig:Garkce_tumors_3D} illustrates the evolution of the tumor and the nutrient concentration in three dimensions. As in the two dimensional case, the tumor grows towards regions with higher nutrient concentrations.

\begin{figure}[htp!]
    \centering
    \includegraphics[width=0.58\linewidth]{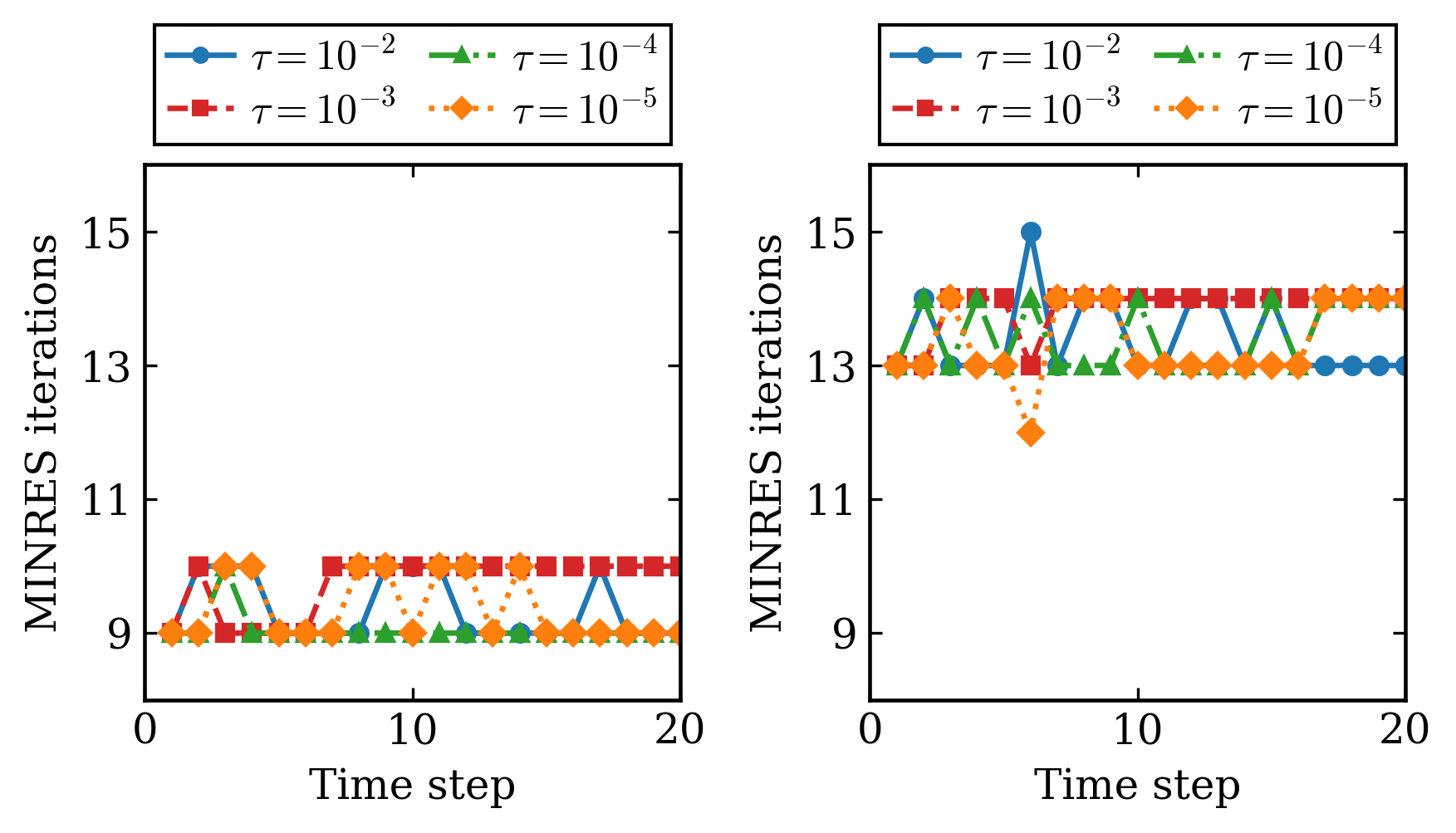}
    \caption{Ex.~\ref{subsec:garcke_smooth} \texttt{MINRES} iterations over $20$ time steps for the solution of system \eqref{GLSS_nonsmooth_matrix_2}, for varying time-step size $\tau$, with initial condition (a). For the solution of system \eqref{GLSS_smooth_matrix}, \texttt{Bi-CGSTAB} is used with \texttt{AMG} preconditioning (left) and with \texttt{ILU} preconditioning (right).}
    \label{fig:GLSS_iter_sigma_3D}
\end{figure}

\subsubsection{Nonsmooth GLSS}\label{subsec:garcke_nonsmooth}
In this section we present numerical results for the nonsmooth variant of the GLSS model. For fixed $c$, the regularized system is smooth and admits a standard Newton linearization, which we solve with \texttt{Bi-CGSTAB} preconditioned by the block preconditioner described in Section~\ref{sec:prec_ns}. We use the same parameters as in Table~\ref{tab:garcke_params}, except that the time-step size is set to $\tau = 10^{-5}$. As $c \to 0$, the condition number of the linearized system grows as $\mathcal{O}(c^{-1})$, since the Moreau--Yosida penalty term introduces a contribution of order $c^{-1}$ to the system matrix. Consequently, a smaller time step is required to maintain Newton convergence at the final regularization level $c_{\max} = 10^{-5}$; with $\tau = 10^{-5}$, the linearized system remains sufficiently well-conditioned for \texttt{Bi-CGSTAB} to converge within a reasonable number of iterations.
 
To obtain an accurate approximation of the obstacle solution at the first time steps, we decrease $c$ through a sequence of levels $ c \in \bigl\{10^{-1},\,10^{-2},\,10^{-3},\,10^{-4},\,10^{-5}\bigr\},$ using the solution at the previous level as the initial guess for the next. In practice, we apply this full continuation only during the first five time steps, after which the solution from the previous time step provides a sufficiently accurate initial guess for the Newton solver at $c = c_{\max} = 10^{-5}$ directly. From step six onwards, only the final level $c = c_{\max}$ is solved at each time step, reducing the computational cost per step significantly.  The initial condition is the four-finger shape (b) defined in Section~\ref{subsec:garcke_smooth}. 

\begin{figure}[htp!]
    \centering
    \includegraphics[width=0.23\linewidth]{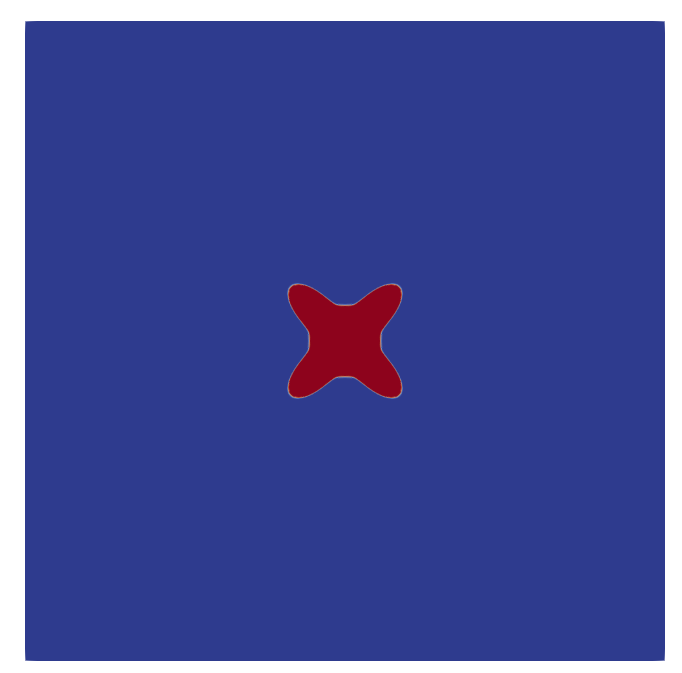}
    \includegraphics[width=0.23\linewidth]{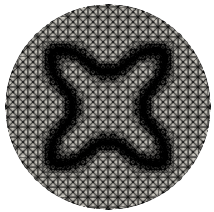}
    \includegraphics[width=0.33\linewidth]{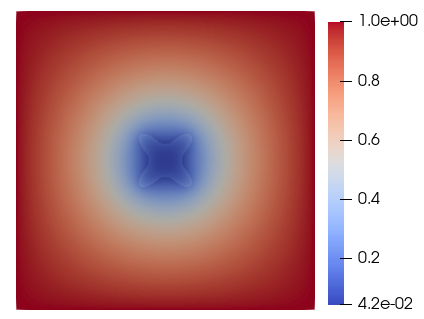}

    \includegraphics[width=0.23\linewidth]{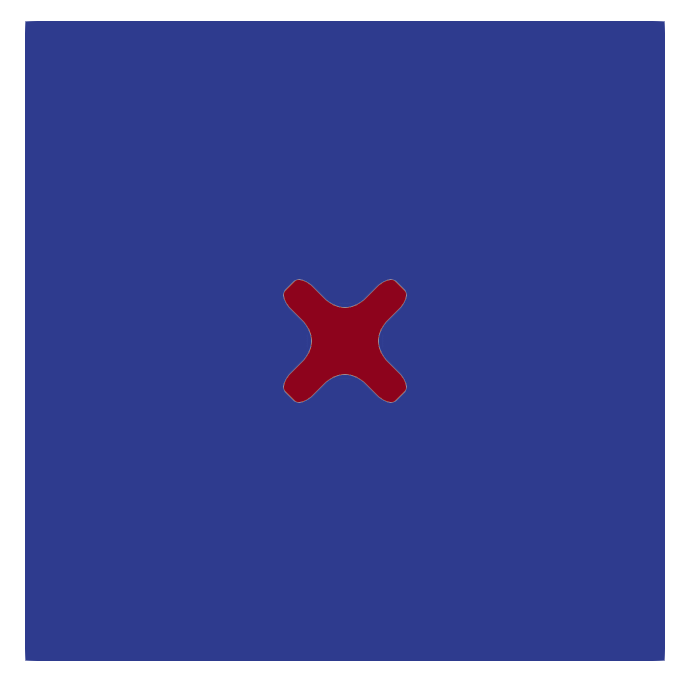}
    \includegraphics[width=0.23\linewidth]{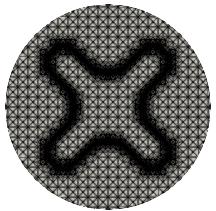}
    \includegraphics[width=0.33\linewidth]{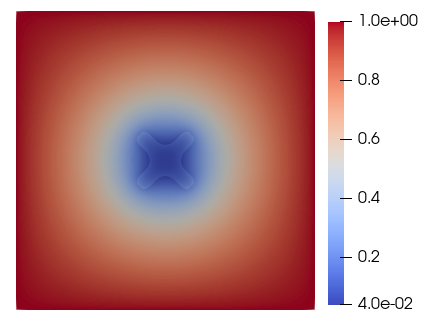}

    \caption{Ex.~\ref{subsec:garcke_nonsmooth} Tumor (left), mesh (middle) and nutrient (left) concentration for time $t = 0.00001$ (top) and  $t = 0.5$ (bottom) with the initial condition (c). }
    \label{fig:GLSS_non_u_sigma_c}
\end{figure}

Figure~\ref{fig:GLSS_non_u_sigma_c} shows the tumor morphology, adaptive mesh, and nutrient concentration at times $t=10^{-5}$ and $t=0.5$ for the nonsmooth GLSS model with the four-finger initial condition~(b). The behavior  of the tumor and nutrient is in close agreement with the smooth formulation. The adaptive mesh correctly tracks the diffuse interface throughout the simulation, concentrating refinement along the tumor boundary while leaving the bulk regions coarse. 

\begin{figure}[htp!]
    \centering
    \includegraphics[width=0.88\linewidth]{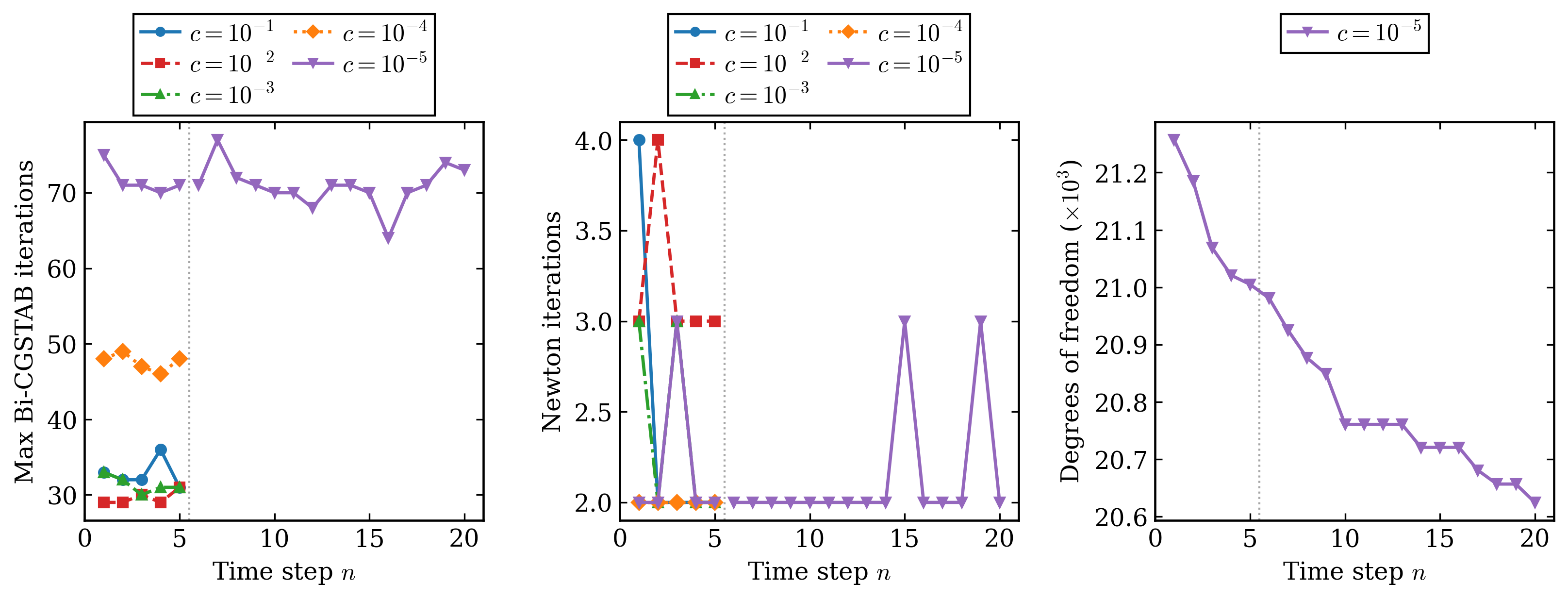}
    \caption{Ex.~\ref{subsec:garcke_nonsmooth} Iteration counts over $20$ time steps for varying time-step size $\tau$: maximum \texttt{Bi-CGSTAB} iterations (left), Newton iterations (center), and degrees of freedom (right) with the initial condition (b).}
    \label{fig:GLSS_noniter_b}
\end{figure}

Figure~\ref{fig:GLSS_noniter_b} presents the maximum \texttt{Bi-CGSTAB} iteration count, the Newton iteration count, and the number of degrees of freedom during the first $20$ time steps for the nonsmooth GLSS model with the four-finger initial condition. As expected, the \texttt{Bi-CGSTAB} solver requires noticeably more iterations at $c=10^{-5}$ than at the intermediate continuation levels, since the linearized systems become increasingly ill-conditioned as $c$ decreases. Despite this, the Newton method converges within at most four iterations throughout the simulation.

Table~\ref{tab:nonsmooth_10k_single} reports the solver statistics over the first $10,000$ time steps. The maximum \texttt{Bi-CGSTAB} iteration count is $251$, while the adaptive mesh contains, on average, approximately $20,173$ degrees of freedom. These results indicate that both the block preconditioner and the adaptive mesh refinement strategy remain stable and efficient over long simulations. Moreover, the phase-field variable remains very close to the obstacle constraints throughout the computation, with $|u_h^n|_{\infty} \leq 1.000054$. The small deviation from unity is caused by the finite Moreau--Yosida regularization parameter $c_{\max}=10^{-5}$. To assess the effectiveness of the proposed preconditioner, we repeated the same experiments without preconditioning. In this case, the \texttt{Bi-CGSTAB} solver failed to converge within the prescribed maximum of $500$ iterations.

Table~\ref{tab:phi_minmax} compares the minimum and maximum values of $u_h^n$ for the smooth and nonsmooth formulations at selected time steps with $\tau=10^{-5}$. In the smooth model, the phase-field exhibits overshoots of up to $|u_h^n|_\infty \approx 1.037$, while the minimum values approach $-1$. In contrast, the nonsmooth formulation keeps $|u_h^n|$ within approximately $5\times10^{-5}$ of the obstacle bounds at all reported time steps, demonstrating the effectiveness of the Moreau--Yosida regularization in enforcing the constraints. The slight violation of the lower bound in the nonsmooth case ($\min u_h^n \approx -1.000041$ at step~$20$) is entirely due to the finite value of the regularization parameter and disappears in the limit $c \rightarrow 0$.

\begin{table}[htbp]
\centering
\caption{Summary over the first $10\,000$ time steps for nonsmooth GLSS model.}
\label{tab:nonsmooth_10k_single}
\small
\setlength{\tabcolsep}{5pt}
\begin{tabular}{cccccc}
\hline
Max \texttt{Bi-CGSTAB} & Max Newton
& $\min\,u_h^n$ & $\max\,u_h^n$ & Avg DOFs & Max DOFs \\
\hline
251 & 16
& $-1.000052$ & $1.000054$ & 20\,173 & 21\,257 \\
\hline
\end{tabular}
\end{table}

\begin{table}[htp!]
\centering
\caption{Minimum and maximum values of $u_h^n$ for smooth and nonsmooth GLSS model, $\tau=10^{-5}$. The nonsmooth model enforces $|u|\leq 1$ via Moreau--Yosida regularization with $c_{\max}=10^{-5}$.}
\label{tab:phi_minmax}
\small
\setlength{\tabcolsep}{4pt}
\begin{tabular}{ll rrrrr}
\hline
Time step & & 1 & 5 & 10 & 15 & 20 \\
\hline
min & smooth    & $-1.000000$ & $-1.000000$ & $-1.000000$ & $-1.000000$ & $-1.000000$ \\
    & nonsmooth & $-1.000000$ & $-1.000000$ & $-1.000027$ & $-1.000038$ & $-1.000041$ \\
\hline
max & smooth    & $1.019581$ & $1.036436$ & $1.036828$ & $1.034482$ & $1.031603$ \\
    & nonsmooth & $1.000049$ & $1.000052$ & $1.000051$ & $1.000051$ & $1.000051$ \\
\hline
\end{tabular}
\end{table}

\section{Conclusions}
\label{sec:concl}
In this work, we have developed and analyzed block Schur complement preconditioners for diffuse interface models of tumor growth: the smooth and nonsmooth variants of the Garcke--Lam--Sitka--Styles (GLSS) model and the Hilhorst--Kampmann--Nguyen--Zee (HKNZ) model. For the smooth GLSS model, the Newton linearization yields a nonsymmetric saddle-point system at each time step, for which we proposed a block-triangular preconditioner based on a Schur complement approximation. Incomplete LU factorizations are used in two dimensions, while algebraic multigrid are employed in three dimensions with \texttt{Bi-CGSTAB}. For the nonsmooth variant, the obstacle potential is regularized via the Moreau--Yosida technique.  For the symmetric HKNZ model, a \texttt{MINRES} solver preconditioned by \texttt{AMG} achieves mesh-independent iteration counts across varying grid sizes and proliferation parameters. Taken together, the numerical results confirm that the proposed preconditioners are robust with respect to model parameters, time-step sizes, and spatial resolution in both two and three dimensions.

\section*{Acknowledgments}
This work was partially supported by the Deutsche Forschungsgemeinschaft (DFG) through Research Unit FOR 5387 POPULAR (Project No.~461909888).

\bibliographystyle{siamplain}
\bibliography{CHtumor_bib}

\end{document}